\documentclass[twocolumn]{autart}

\usepackage{graphicx,amsmath,amsfonts,graphics,subfigure}          % Include this line if your
\usepackage{indentfirst}
\usepackage{xcolor}
\usepackage{epstopdf}
\usepackage[authoryear,longnamesfirst]{natbib}
\usepackage{hyperref}
\usepackage[none]{hyphenat}

\begin{document}

\begin{frontmatter}

\title{Stochastic Bounded Real Lemma and $H_{\infty}$ Control \\ of Difference Systems in Hilbert Spaces}
\thanks[footnoteinfo]{Corresponding author.}

\author[sdnu]{Cheng'ao Li}\ead{$lichengao103@163.com$},
\author[sdnu]{Ting Hou \thanksref{footnoteinfo}}\ead{$ht\_math@sina.com$},
\author[sdust]{Weihai Zhang}\ead{$w\_hzhang@163.com$},
\author[scut]{Feiqi Deng}\ead{$aufqdeng@scut.edu.cn$},

\address[sdnu]{School of Mathematics and Statistics, Shandong Normal University,
Jinan 250014, Shandong Province, China}
\address[sdust]{College of Electrical Engineering and Automation, Shandong University of Science and Technology,\\
Qingdao 266590, Shandong Province, China}
\address[scut]{School of Automation Science and Engineering, South China University of Technology,\\
Guangzhou 510640, Guangdong Province, China}


\begin{abstract}
This paper mainly establishes the finite-horizon stochastic bounded real lemma, and then solves the $H_{\infty}$ control problem
for discrete-time stochastic linear systems defined on the separable Hilbert spaces,
thereby unifying the relevant theoretical results previously confined to the Euclidean space $\mathbb{R}^n$.
To achieve these goals, the indefinite linear quadratic (LQ)-optimal control problem is firstly discussed.
By employing the bounded linear operator theory and the inner product,
a sufficient and necessary condition for the existence of
a linear state feedback LQ-optimal control law is derived,
which is closely linked with the solvability of the backward Riccati operator equation with a sign condition.
Based on this, stochastic bounded real lemma is set up to facilitate the $H_{\infty}$ performance of the disturbed system in Hilbert spaces.
Furthermore, the Nash equilibrium problem associated with two parameterized quadratic performance indices is worked out,
which enables a uniform treatment of the $H_{\infty}$ and $H_2/H_{\infty}$ control designs by selecting specific values for the parameters.
Several examples are supplied to illustrate the effectiveness of the obtained results,
especially the practical significance in engineering applications.
\end{abstract}

\end{frontmatter}

\section{Introduction}

Optimal control of stochastic systems has emerged as a critical area of research
due to its extensive applications across various fields including engineering, finance, economics, and so on.
The stochastic nature of this kind of systems, compound by uncertainties and noise,
presents significant challenges in ensuring stability and achieving the desired performance.
Conventional finite-dimensional models, often represented in $\mathbb{R}^n$,
have been extensively researched in optimal control over the past decades
(see, \citet{2004Stochastic, dragan2006mathematical, 2021Stability, 2008Generalized}).
However, these models exhibit inherent limitations when dealing with the complexities of real-world systems.
As is well-known, some systems, such as delay systems, distributed parameter systems or partial differential systems and so on,
cannot be adequately modeled within a finite-dimensional framework.
For instance, fluid systems, structural vibration systems, and quantum mechanical systems,
which have the temporal and spatial variation characteristics,
are typically described by the infinite-dimensional dynamics.
These systems require a more sophisticated mathematical framework, such as Hilbert spaces, to accurately model their states and dynamics
(see, \citet{chueshov2002introduction, dutailly2014hilbert, 2003Infinite, Fedorov2003Linear}).
Since the finite-dimensional models may oversimplify the complexities at hand,
which will result in less-than-optimal control strategies,
they often cannot remain effective across a range of scenarios.
Therefore, it is important to recognize that in certain cases, these models may not be sufficient,
and more advanced techniques are necessary to find the optimal control.

This paper aims to set up the bounded real lemma and tackle the $H_{\infty}$ control problem over the finite-horizon
for discrete-time stochastic systems within the scope of Hilbert spaces.
Our findings are consistent with the well-established principles of the corresponding finite-dimensional control theories,
while overcoming the challenges encountered when modeling systems with high complexity.
Considering that the $H_{\infty}$ control and the indefinite linear quadratic (LQ)-optimal control are two closely related issues,
our foremost objective is to form a relatively complete theory on the finite-horizon indefinite LQ-optimal control
within the framework of Hilbert spaces.
To achieve this, we follow a similar approach used in finite-dimensional studies of stochastic LQ-optimal control problems
(see, \citet{chen1998stochastic, chen2000stochastic, rami2000linear, wang2024stochastic}),
where the LQ-optimal control problem is typically transformed into the solvability of a class of Riccati equations.
By employing the linear operator theory in \citet{akhiezer2013theory, da2014stochastic, ungureanu2004uniform, ungureanu2013stability},
the matrices included in finite-dimensional systems are replaced by bounded linear operators,
which allows us to substitute the traditional backward Riccati equation with the backward Riccati operator equation.
Since the matrices can be considered as a special case of bounded linear operators,
our adoption of the backward Riccati operator equation not only preserves the essential theoretical results
established in the finite-dimensional context,
but also extends their applicability to the infinite-dimensional setting,
thereby providing a unified framework for addressing the LQ-optimal control problem.

Moreover, also within the Hilbert spaces framework, stochastic bounded real lemma,
which plays a crucial role in characterizing the $H_{\infty}$ performance criteria
(see, \citet{2010Mathematical, hinrichsen1998stochastic, XIAO2024111827}),
is built upon the developed indefinite LQ-optimal control theories.
In addition to this, the discussion of optimal control strategies of a robust system is formulated as a Nash game problem.
A sufficient and necessary condition for the existence of Nash equilibria is proposed
by conducting a thorough analysis on the solvability of the coupled backward Riccati operator equations.
All of these analyses not only deepen our understanding of stochastic control theories in the infinite-dimensional setting,
but also provide new insights for the design of robust control strategies for complex systems.

The main contributions of this paper are summarized as follows:
\begin{itemize}
\item  We generalize the classical finite-horizon LQ-optimal control, $H_{\infty}$ analysis, and the Nash game theory to Hilbert spaces.
       As Hilbert space serves as a natural framework for representing infinite-dimensional systems,
       our research constructs a unified theoretical framework that is applicable to both finite-dimensional and infinite-dimensional cases.
\item  To better analyze the $H_{\infty}$ performance of difference systems with disturbance in an infinite-dimensional context,
       we establish the finite-horizon stochastic bounded real lemma based on the bounded linear operator theory.
       To the best of our knowledge, this is the first time that a sufficient and necessary condition is presented
       within the framework of Hilbert spaces to ensure that the $H_{\infty}$ norm of the associated perturbation operator is below a prescribed level $\gamma>0$.
\item  Two practical examples (Examples \ref{example-1} and \ref{example-2}) are included to confirm the significant engineering application value
       of the control theories developed for the system whose state space is an infinite-dimensional space.
       One involves applying convolution operators to process and control audio signals,
       and the other is concerned with designing an optimal control strategy for a heat transfer system
       such that the temperature deviation is minimized while the energy usage is restricted effectively to avoid unnecessary power consumption.
\end{itemize}

The remainder is organized as follows:
Section 2 is with the formulation of the difference system model in Hilbert spaces.
In Section 3, the well-posedness of the indefinite LQ-optimal control problem is first addressed.
Then, a sufficient and necessary condition is presented for the existence of a linear state feedback LQ-optimal law.
Based on this, the finite-horizon stochastic bounded real lemma is built in Section 4.
Section 5 explores the Nash game problem associated with two parameterized quadratic performance indices,
and thus resolves $H_{\infty}$ and $H_2/H_{\infty}$ control problems uniformly.
Section 6 concludes this paper with a brief summary.

\section{Notations and the Model Formulation}

Let $H$, $U$, $V$ and $Z$ be real separable Hilbert spaces
and $L(H,U)$ be the Banach space of all bounded linear operators transforming $H$ into $U$.
$\varepsilon(H)$ is the Banach subspace of $L(H): = L(H,H)$, formed by all self-adjoint operators.
Denote the inner product and norm in Hilbert space $H$ as $\langle \cdotp, \cdotp \rangle_H$ and $\Arrowvert \cdotp\Arrowvert_H$, respectively.
As usual, $*$ means either the adjoint of a bounded linear operator or the dual of a Banach space.

$A\in L(H)$ is said to be nonnegative, and denoted by $A\geq 0$,
if $A$ is self-adjoint and $\langle Ax, x \rangle_H \geq 0$ for any $x\in H$.
$A\in L(H)$ is said to be positive, and denoted by $A>0$, if there exists $\vartheta >0$ such that $A-\vartheta I_{H} \geq 0$,
where $I_{H}$ is the identity operator on $H$.
It is easy to verify that a positive bounded linear operator is an invertible operator.
Set $(\Omega, \mathcal{F}, P)$ be a probability space.
$\xi$ is a real or $H$-valued random variable on $\Omega$ and its mean value (expectation) is expressed by $E[\xi]$.
Let $\mathcal{Q}=\{0,1,\cdots,N\}$, and $\mathcal{Q}^{+}=\{0,1,\cdots,N+1\}$.
$K$ is assumed to be a real Banach space with the norm $\Arrowvert \cdotp \Arrowvert_{K}$,
and $K^{\mathcal{Q}}$ denotes the set of the sequences $q=\{q(i) | q(i)\in K\}_{i\in \mathcal{Q}}$.
Let $L^2_{K^{\mathcal{Q}}}=\{q\in K^{\mathcal{Q}},\ \Arrowvert q \Arrowvert_{L^2_{K^{\mathcal{Q}}}}=(\sum_{i\in \mathcal{Q}} E[\Arrowvert q(i) \Arrowvert^2_{K}])^{\frac{1}{2}}<+\infty \}$ be a real Banach space with the norm $\Arrowvert \cdotp \Arrowvert_{L^2_{K^{\mathcal{Q}}}}$.
\begin{rem}\label{remake-1}
According to Proposition 3 of \citet{dutailly2014hilbert},
any real separable Hilbert space can be endowed with the structure of a complex separable Hilbert space;
therefore, without loss of generality, the discussion in this paper is restricted to real separable Hilbert spaces.
\end{rem}

On the probability space $(\Omega, \mathcal{F}, P)$, consider the following discrete-time linear system with multiplicative noise:
\begin{equation}\label{1}
	\begin{split}
		 x(k+1)=&A(k)x(k)+B(k)u(k)+(C(k)x(k)\\
			    &+D(k)u(k))\omega(k),\\
		x(0)=&x_0\in H,
	\end{split}
\end{equation}
where $x(k)\in H$ and  $u(k)\in U$ ($k\in \mathcal{Q}$) are respectively the state and control input,
$A(k)$, $C(k)\in L(H)$ and $B(k)$, $D(k)\in L(U,H)$.
The initial value $x_0$ is assumed to be a deterministic vector.
$\{\omega(k),\ k\in \mathcal{Q}\}$ is a sequence of real random variables with
$E[\omega(k)]=0$ and $E[\omega(k_1)\omega(k_2)]=\delta(k_1,k_2)$, $k_1,\ k_2\in \mathcal{Q}$, in which $\delta$ is the Kronecker function.
Denote $\mathcal{F}_k$  the $\sigma$-field generated by $\{\omega(0),\cdots,\omega(k-1)\}$ and $\mathcal{F}_0=\sigma(x_0)$.
\begin{rem}	\label{remake-2}
Stochastic difference equation \eqref{1} describes a type of systems with state and control input dependent noise.
Note that for any $k\in \mathcal{Q}$, $\omega(k)$ is a real random variable.
Then, from \eqref{1}, it can be known that $x(k)$ is a random variable on $(\Omega, \mathcal{F}, P)$ with values in $H$.
\end{rem}
\begin{rem}	\label{remake-3}
Different from the most existing literature dealing with a finite-dimensional state space $\mathbb{R}^n$,
this paper replaces the Euclidean space $\mathbb{R}^n$ with a (potentially infinite-dimensional) real separable Hilbert space.
Correspondingly, matrices $A$, $B$, $C$ and $D$ of suitable dimensions are substituted with bounded linear operators,
still denoted as $A$, $B$, $C$ and $D$ for convenience.
This extension enables the model to handle infinite-dimensional states and control variables,
and thereby covers a wider range of application scenarios.
\end{rem}

\section{Indefinite LQ-optimal Control}

Set $\mathcal{U}$ be the class of controls $u=(u(0),u(1),\cdots ,\\ u(N))$
with $u(k)\in U$ and $E[\Arrowvert u(k)\Arrowvert_U^2]<+\infty$ being $\mathcal{F}_k$-measurable for each $k\in \mathcal{Q}$.
The linear quadratic cost functional associated with \eqref{1} is defined as follows:
\begin{equation}\label{2}
	\begin{split}
J(x_0,u) &=E\bigg[\sum_{k=0}^{N} \Big( \langle M(k)x(k),x(k) \rangle_H\\
&\quad\quad + 2\langle L(k)x(k),u(k)\rangle_U + \langle R(k)u(k),u(k)\rangle_U\Big) \\
&\quad\quad + \langle S(N+1)x(N+1),x(N+1)\rangle_H \bigg],
	\end{split}
\end{equation}
where $x(\cdot)$ is the solution to system \eqref{1} corresponding to the control input $u(\cdot)$.
Moreover, $L(k)\in L(H,U)$, $M(k),\ S(N+1)\in \varepsilon(H)$, $R(k)\in \varepsilon(U)$ for each $k\in \mathcal{Q}$.

The finite-horizon indefinite LQ-optimal control problem is to find an optimal control $u^* \in \mathcal{U}$
such that the linear quadratic cost functional \eqref{2} is minimized
according to system \eqref{1} with the initial $H$-valued vector $x_0$,
and the optimal cost value is represented as $\mathcal{J}(x_0)$, i.e.
\begin{equation}\label{3}
	\begin{split}
		\mathcal{J}(x_0)=\inf_{u \in \mathcal{U}}J(x_0,u)=J(x_0,u^*).
	\end{split}
\end{equation}
\begin{defn}\label{definition-1}
The LQ-optimal control problem \eqref{3} is called to be well-posed, if
$\mathcal{J}(x_0)>-\infty$ for any initial value $x_0 \in H$.
\end{defn}

In the following, we attempt to provide sufficient conditions to ensure the well-posedness
of the indefinite LQ-optimal control problem \eqref{3},
which is not an ordinary work within infinite-dimensional spaces.
To this end, we first introduce the associated backward Riccati operator equation.

Let $Dom(\Pi_k)=\{X\in \varepsilon(H)|R(k)+B(k)^*XB(k)+D(k)^* XD(k) \ is \ invertible \ with \ a \ bounded \ inverse\}$
and define the Riccati operator $\Pi_k:\ Dom(\Pi_k)\mapsto \varepsilon(H)$, $k\in \mathcal{Q}$ as follows:
\begin{equation}\label{4}
	\begin{split}
		\Pi_k(X)=&A(k)^*XA(k)+C(k)^*XC(k)+M(k)\\
		&-[L(k)^*+A(k)^*XB(k)+C(k)^*XD(k)]\\
		&\cdot [R(k)+B(k)^*XB(k)+D^*(k) XD(k)]^{-1}\\
		&\cdot [L(k)+B(k)^*XA(k)+D(k)^*XC(k)].
	\end{split}
\end{equation}
For $k\in \mathcal{Q}$, consider the following backward Riccati operator equation:
\begin{equation}\label{5}
	\left\{\begin{aligned}
		&P(k)=\Pi_k(P(k+1)),\\
		&P(N+1)=S(N+1)\in \ Dom(\Pi_N).
	\end{aligned}
	\right.
\end{equation}
\begin{defn}\label{definition-2}
$\{P(k),\ k\in\mathcal{Q}^{+}\}$ is called a solution to the backward Riccati operator equation \eqref{5},
if $P(k+1)\in Dom(\Pi_k)$ for any $k\in\mathcal{Q}$ and \eqref{5} is satisfied.
\end{defn}
\begin{lem}\label{lemma-1}
Suppose that $\{P(k),\ k\in\mathcal{Q}^{+}\}$ is a solution to the backward Riccati operator equation \eqref{5},
then for any $x_0\in H$ and $u \in \mathcal{U}$, we have that
	\begin{equation*}
		\begin{split}
			J(x_0,u)=&\langle P(0) x_0,x_0\rangle_H+E\bigg[\sum_{k=0}^{N}\langle \mathcal{R}(k)\big(u(k)+\mathcal{R}(k)^{-1}\\
                     & \cdot \mathcal{G}(k) x(k)\big), u(k)+\mathcal{R}(k)^{-1}\mathcal{G}(k) x(k)\rangle_U\bigg],
		\end{split}
	\end{equation*}
where $x(\cdot)$ is the corresponding solution to system \eqref{1}, and for any $k\in \mathcal{Q}$,
$\mathcal{R}(k) = R(k)+B(k)^*P(k+1)B(k)+D(k)^* P(k+1)D(k),$
$\mathcal{G}(k) = L(k)+B(k)^*P(k+1)A(k)+D(k)^*P(k+1)C(k).$
\end{lem}
\begin{pf}
In view of the Riccati operator $\Pi_k$
and the fact that $\{P(k),\ k\in\mathcal{Q}^{+}\}$ is a solution to the backward Riccati operator equation \eqref{5},
one can infer that
\begin{equation}\label{6}
		\begin{split}
			&J(x_0,u)\\
            &=J(x_0,u)+E\bigg[\sum_{k=0}^{N}\Big(\langle P(k+1)x(k+1),x(k+1)\rangle_H \\
			&\ \ \ \ -\langle P(k)x(k),x(k)\rangle_H\Big)+\langle P(0) x_0,x_0\rangle_H\\
			&\ \ \ \ -\langle P(N+1)x(N+1),x(N+1)\rangle_H \bigg]\\
			&=E\bigg[\langle P(0) x_0,x_0\rangle_H+\sum_{k=0}^{N} \Big(\langle \big[M(k)+A(k)^*P(k+1)A(k)\\
			&\ \ \ \ +C(k)^*P(k+1)C(k)\big]x(k),x(k)\rangle_H+2\langle \big[L(k)\\
			&\ \ \ \ +B(k)^*P(k+1)A(k)+D(k)^*P(k+1)C(k)\big]x(k),\\
			&\ \ \ \ \ u(k)\rangle_U+\langle \big[R(k)+B(k)^*P(k+1)B(k)\\
			&\ \ \ \ +D(k)^*P(k+1)D(k)\big]u(k),u(k)\rangle_U\\
            &\ \ \ \ -\langle P(k)x(k),x(k)\rangle_H\Big)\bigg].
        \end{split}
\end{equation}
Further, considering the definitions of $\mathcal{R}(k)$ and $\mathcal{G}(k)$, it follows from \eqref{6} that
\begin{equation*}
		\begin{split}
			&J(x_0,u)\\
            &=\langle P(0) x_0,x_0\rangle_H+E\bigg[\sum_{k=0}^{N} \langle \mathcal{G}(k)^*\mathcal{R}(k)^{-1} \mathcal{G}(k)x(k),\\
			&\ \ \ \ \ x(k)\rangle_H+2\langle \mathcal{G}(k) x(k),u(k)\rangle_U+\langle \mathcal{R}(k)u(k),u(k)\rangle_U\bigg]\\
            &=\langle P(0) x_0,x_0\rangle_H+E\bigg[\sum_{k=0}^{N}\langle \mathcal{R}(k)\big(u(k)+\mathcal{R}(k)^{-1}\\
            &\ \ \ \ \ \cdot \mathcal{G}(k) x(k)\big), u(k)+\mathcal{R}(k)^{-1}\mathcal{G}(k) x(k)\rangle_U\bigg].
        \end{split}
\end{equation*}
Hence, Lemma \ref{lemma-1} is proved.
\end{pf}

Denote $\ell^2=\{(\xi_1,\xi_2,\cdots)|\sum_{i=1}^{+\infty}|\xi_i|^2<+\infty\}$.
It is well-known that any one separable Hilbert space $H$ is isometric isomorphic to $\ell^2$
(for example, see Theorem 5.4 in \citet{conway1990}).
Thus, for any $x\in H$, one has that $\Arrowvert x \Arrowvert_{H}=\Arrowvert Tx \Arrowvert_{\ell^2}<+\infty$,
in which $T$ is the isometric isomorphic mapping from $H$ to $\ell^2$.
Attributing to Lemma \ref{lemma-1} and the finiteness of the norm of vectors in separable Hilbert spaces,
the following proposition is drawn:
\begin{prop}\label{proposition-1-2}
If the backward Riccati operator equation \eqref{5} has a solution $\{P(k),\ k\in\mathcal{Q}^{+}\}$
such that the operator family $\{\mathcal{R}(k),\ k\in\mathcal{Q}\}$ is uniformly positive,
i.e., for any $k\in\mathcal{Q}$, $\mathcal{R}(k)>0$,
then the LQ-optimal control problem \eqref{3} is well-posed.
\end{prop}
\begin{pf}
Under the assumption that $\{\mathcal{R}(k),\ k\in\mathcal{Q}\}$ is uniformly positive,
from Lemma \ref{lemma-1}, it follows that $J(x_0,u)\geq \langle P(0) x_0,x_0\rangle_H$ for all $u \in \mathcal{U}$.
Hence, $\mathcal{J}(x_0)\geq \langle P(0) x_0,x_0\rangle_H>-\infty$.
This ends the proof.
\end{pf}

If the optimal solution is specified to be searched in a state feedback form,
the following theorem provides a sufficient and necessary condition for the solvability of
the LQ-optimal control problem \eqref{3} associated with system \eqref{1}.
\begin{thm}\label{theorem-1}
For system \eqref{1}, the LQ-optimal control problem \eqref{3} admits a unique optimal control
and that control is in a linear state feedback form
if and only if (iff) the backward Riccati operator equation \eqref{5} has a solution $\{P(k),\ k\in\mathcal{Q}^{+}\}$
such that the operator family $\{\mathcal{R}(k),\ k\in\mathcal{Q}\}$ is uniformly positive.
Moreover, the optimal control is explicitly given by $u^*(k)=-\mathcal{R}(k)^{-1}\mathcal{G}(k) x(k)$, $k\in\mathcal{Q}$,
and the optimal cost value $\mathcal{J}(x_0)=J(x_0,u^*)=\langle P(0) x_0,x_0\rangle_H$.
\end{thm}
\begin{pf}
Sufficiency.
Assume that the backward Riccati operator equation \eqref{5} has a solution $\{P(k),\ k\in\mathcal{Q}^{+}\}$
such that the operator family $\{\mathcal{R}(k),\ k\in\mathcal{Q}\}$ is uniformly positive,
and we intend to search the optimal solution in a state feedback form.
First of all, from Proposition \ref{proposition-1-2},
it yields that the LQ-optimal control problem \eqref{3} is well-posed
and $\mathcal{J}(x_0)\geq \langle P(0) x_0,x_0\rangle_H$.
On the other hand,
by constructing a linear state feedback control $u^*(k)\!=-\mathcal{R}(k)^{-1} \mathcal{G}(k) x(k)\ (k\in \mathcal{Q})$
and applying Lemma \ref{lemma-1},
we have that $J(x_0, u^*)=\langle P(0) x_0, x_0 \rangle_H$,
and thereby $\mathcal{J}(x_0)\leq\langle P(0) x_0, x_0\rangle_H$.
Consequently, $\mathcal{J}(x_0)=\langle P(0) x_0, x_0 \rangle_H$ and $u^*(k)$ is an optimal control.
Suppose that the LQ-optimal control problem \eqref{3} admits another linear state feedback control $\hat{u}(k)\ (k\in\mathcal{Q})$
such that $\mathcal{J}(x_0)=J(x_0,\hat{u})$ holds.
From Proposition  \ref{proposition-1-2}, one knows that $\mathcal{J}(x_0)>-\infty$,
and $J(x_0,u^*)-J(x_0,\hat{u})=\mathcal{J}(x_0)-\mathcal{J}(x_0)=0$ is then followed.
By Lemma \ref{lemma-1}, it concludes that $\hat{u}(k)=u^*(k)$, $k\in\mathcal{Q}$.
	
Necessity.
Assume that the LQ-optimal control problem \eqref{3} associated with system \eqref{1}
admits a unique optimal control and that control is in a linear state feedback form.
Firstly, we prove the existence of the solution to the backward Riccati operator equation \eqref{5} by contradiction.
If it is false, then there exists $k_0 \in \mathcal{Q}$ such that
$P(k_0+1)\notin Dom(\Pi_{k_0})$ and $P(k+1)\in Dom(\Pi_{k})$ for any $k>k_0$,
i.e., there exists $\widetilde{u}(k_0) \in U$, $\widetilde{u}(k_0)\neq0$
satisfying $\langle \mathcal{R}(k_0)\widetilde{u}(k_0),\widetilde{u}(k_0)\rangle_U=0$.
Take $x_0=0$, and suppose that $u^0=\{u^0(k),\ k\in\mathcal{Q}\}$ is the linear state feedback control law
solves the corresponding LQ-optimal control problem \eqref{3}, i.e. $J(0,u^0)=\mathcal{J}(0)=0$.
Now, construct a nonlinear feedback control sequence $u^\Upsilon=\{u^\Upsilon(k)$, $k\in\mathcal{Q}\}$ as follows:
	\begin{equation*}
		    u^\Upsilon(k)=\left\{\begin{aligned}
			&0,                   &k<k_0, \\
			&\widetilde{u}(k_0),  &k=k_0,\\
			&-\mathcal{R}(k)^{-1}\mathcal{G}(k)x(k), &k>k_0.
		\end{aligned}\right.
	\end{equation*}
Note that for every $k>k_0$,
the condition $P(k+1)\in Dom(\Pi_{k})$ ensures the existence of $\mathcal{R}(k)^{-1}$.
Then, $u^\Upsilon$ is well defined, and $u^\Upsilon\in \mathcal{U}$.
Hence, for system \eqref{1} with $x_0=0$, under this control input $u^\Upsilon$,
the state satisfies $x(k)=0$ for all $k \leq k_0$.
Furthermore, by a derivation similar to that of \eqref{6} on the part of the trajectory ($k>k_0$)
where the solution to the backward Riccati operator equation \eqref{5} exists,
it can be deduced that
\begin{equation*}
		\begin{split}
			&J(0,u^\Upsilon)\\
            =&E\bigg[\langle R(k_0)\widetilde{u}(k_0),\widetilde{u}(k_0)\rangle_U+\sum_{k=k_0+1}^{N} \Big( \langle M(k)x(k),x(k) \rangle_H\\
            &\ \ \ \ + 2\langle L(k)x(k),u^\Upsilon(k)\rangle_U+ \langle R(k)u^\Upsilon(k),u^\Upsilon(k)\rangle_U \Big)\\
            &\ \ \ \ +\langle S(N+1)x(N+1),x(N+1)\rangle_H  \bigg]\\
            =&\langle \mathcal{R}(k_0)\widetilde{u}(k_0),\widetilde{u}(k_0)\rangle_U+E\bigg[\sum_{k=k_0+1}^{N}\langle \mathcal{R}(k)\big(u^\Upsilon(k)\\
            &+\mathcal{R}(k)^{-1}\mathcal{G}(k) x(k)\big), u^\Upsilon(k)+\mathcal{R}(k)^{-1}\mathcal{G}(k) x(k)\rangle_U\bigg]\\
            =&\langle \mathcal{R}(k_0)\widetilde{u}(k_0),\widetilde{u}(k_0)\rangle_U\\
            =&0,
		\end{split}
	\end{equation*}
which contradicts with the assumption that $\mathcal{J}(0)$ is achieved via the linear state feedback control law $u^0$.
Therefore, the backward Riccati operator equation \eqref{5} has a solution $\{P(k),\ k\in\mathcal{Q}^{+}\}$.

Next, we will show that the operator family $\{\mathcal{R}(k),\ k\in\mathcal{Q}\}$ is uniformly positive.
Notice that for any $k\in\mathcal{Q}$, $\mathcal{R}(k)^{-1}$ exists, that is,
there is no nonzero $u(k)\in U$ ($k\in\mathcal{Q}$) satisfying $\mathcal{R}(k) u(k)=0$.
If the operator family $\{\mathcal{R}(k),\ k\in\mathcal{Q}\}$ is not uniformly positive,
then there exists $k_1 \in \mathcal{Q} $, $\check{u}(k_1) \in U$ and $\check{u}(k_1)\neq0$,
such that $\langle \mathcal{R}(k_1)\check{u}(k_1),\check{u}(k_1)\rangle_U<0$.
Following a similar argument as above, a contradiction will be led to.
So the remaining part of the proof is omitted.
\end{pf}
\begin{rem}\label{remake-4}
Since the matrix operators can also be viewed as bounded linear operators,
Theorem \ref{theorem-1} extends the findings of the finite-horizon LQ-optimal control problem \eqref{3}
defined within finite-dimensional spaces.
\end{rem}

From the discussion above, it is clear that there is a close connection
between the backward Riccati operator equation \eqref{5} and the finite-horizon indefinite LQ-optimal control problem.
Below, the Schur complement of the bounded linear operators is employed
to take a deeper look into the properties of the solution to the backward Riccati operator equation \eqref{5}.

Let $H$ and $U$ be real separable Hilbert spaces, and $M_{11}\in \varepsilon(H)$, $M_{22}\in \varepsilon(U)$ and $M_{21}\in L(H,U)$.
It is not difficult to prove that the product space $H \times U$ is a Hilbert space with the inner product $\langle (x_0,u_0),(x_1,u_1)\rangle_{H \times U} =\langle x_0,x_1\rangle_{H}+\langle u_0,u_1\rangle_{U}$, $\forall (x_0,u_0),\\ (x_1,u_1) \in H \times U$.
For any $(x,u) \in H \times U$, define the operator $M$ as follows:
\begin{equation*}
	\begin{split}
		 &M=\left( \begin{matrix}
			M_{11} & M_{21}^{*}\\ M_{21} &M_{22}
		\end{matrix} \right) :H \times U\mapsto H \times U, \\
		&M(x,u)=(M_{11}x+M_{21}^{*}u, M_{21}x+ M_{22}u).
	\end{split}
\end{equation*}
Obviously, $M \in\varepsilon (H \times U)$.
\begin{defn}\label{definition-3}(\cite{schurgc2002})
If $M_{22}>0$, then the bounded linear operator $M | M_{22}=M_{11}-M_{21}^{*}M_{22}^{-1}M_{21}$ is well defined
and called the Schur complement of $M_{22}$ in $M$.
\end{defn}

The following lemma shows that the standard Schur's Complement Lemma
(for example, see \citet{boyd1994}, page 7)
remains valid for the bounded linear operators.
\begin{lem}\label{lemma-2}(\cite{2013Global})
Assume that $M$ is a bounded linear operator with $M_{22}>0$, then $M>0$ $( M\geq 0)$ iff
$M | M_{22}>0$ $(M | M_{22}\geq 0)$.
\end{lem}

For $k\in \mathcal{Q}$, set $\varPsi_k=\left( \begin{matrix}M(k) & L(k)^{*}\\ L(k) &R(k)\end{matrix} \right)$.
By using Lemma \ref{lemma-2},
we can derive the following existence condition of the solution to the backward Riccati operator equation \eqref{5}.
\begin{thm}\label{theorem-2}
If $S(N+1)\geq 0$, $R(k)>0$, and $\varPsi_k\geq 0$ hold for any $k\in\mathcal{Q}$,
then the backward Riccati operator equation \eqref{5} admits a uniformly nonnegative solution $\{P(k),\ k\in\mathcal{Q}^{+}\}$.
\end{thm}
\begin{pf}
In view of $P(N+1)=S(N+1)\geq 0$ and $R(N)>0$, one has that $\mathcal{R}(N)>0$, so $P(N+1)\in Dom(\Pi_N)$.
From Lemma \ref{lemma-2}, it is known that $P(N)=\Pi_N(P(N+1)) \geq 0$ iff
    \begin{equation*}
		\Xi_N=\left( \begin{matrix}
			\Xi_N^{1} & {\Xi_N^{2}}^*\\ \Xi_N^{2} & \Xi_N^{3}
		      \end{matrix} \right)\geq 0,
	\end{equation*}
where
$\Xi_N^{1}=[A(N)^*P(N+1)A(N)+C(N)^*P(N+1)C(N)+M(N)],$
$\Xi_N^{2}=[B(N)^*P(N+1)A(N)+D(N)^*P(N+1)C(N)+L(N)],$
and $\Xi_N^{3}=[B(N)^*P(N+1)B(N)+D(N)^*P(N+1)D(N)+R(N)].$
Suppose that $\Xi_N$ is not nonnegative, then there is $(x_N,u_N) \in H \times U$
such that $\langle \Xi_N(x_N,u_N) , (x_N,u_N)\rangle <0,$
that is,
	\begin{equation}\label{eqa}
		\begin{split}
			&\langle \Xi_N(x_N,u_N) , (x_N,u_N)\rangle_{H \times U}\\
            &=\langle (\Xi_N^{1}x_N+{\Xi_N^{2}}^*u_N , \Xi_N^{2}x_N+\Xi_N^{3}u_N),(x_N,u_N) \rangle_{H \times U} \\
			&=\langle M(N)x_N,x_N\rangle_H+\langle L(N)^{*}u_N,x_N\rangle_H\\
			&\ \ \ \ +\langle L(N)x_N,u_N\rangle_U+\langle R(N)u_N,u_N\rangle_U\\
			&\ \ \ \ +\langle R(N)u_N,u_N\rangle_U+\langle P(N+1)[A(N)x_N+B(N)u_N],\\
			&\ \ \ \ \ [A(N)x_N+B(N)u_N]\rangle_H+\langle P(N+1)[C(N)x_N\\
			&\ \ \ \ \ +D(N)u_N], [C(N)x_N+D(N)u_N] \rangle_H \\
			&<0.
		\end{split}
	\end{equation}
Furthermore, directly from $\varPsi_N=\left( \begin{matrix}
		M(N) & L(N)^{*}\\ L(N) &R(N)
	\end{matrix} \right) \geq 0$,
we come to the conclusion that
	\begin{equation}\label{eqb}
		\begin{split}
			&\langle M(N)x_N,x_N\rangle_H+\langle L(N)^{*}u_N,x_N\rangle_H\\
			&+\langle L(N)x_N,u_N\rangle_U+\langle R(N)u_N,u_N\rangle_U \geq 0.
		\end{split}
	\end{equation}
Combining \eqref{eqa} with \eqref{eqb}, one has that
$\langle P(N+1)[A(N)x_N+B(N)u_N],[A(N)x_N+B(N)u_N]\rangle_H+\langle P(N+1)[C(N)x_N+D(N)u_N],[C(N)x_N+D(N)u_N]\rangle_H<0,$
which contradicts with $P(N+1) \geq 0$.
Thus, $\Xi_N \geq 0$ and  $P(N)=\Pi_N(P(N+1)) \geq 0$ by Lemma \ref{lemma-2}.
	
Similarly, we can conclude that $P(N)\in Dom(\Pi_{N-1})$ and $P(N-1) \geq 0$.
By inductive method, it follows that $P(k+1)\in Dom(\Pi_k)$ for any $k\in\mathcal{Q}$
and $\{P(k),\ k\in\mathcal{Q}^{+}\}$ is uniformly nonnegative.
The proof of Theorem \ref{theorem-2} is complete.
\end{pf}

It is noted that Theorem \ref{theorem-2} ensures not only the existence and uniformly nonnegativity of
the solution to the backward Riccati operator equation \eqref{5},
but also the uniformly positivity of the operator family $\{\mathcal{R}(k),\ k\in\mathcal{Q}\}$.
Therefore, from Theorems \ref{theorem-1} and \ref{theorem-2}, the following corollary can be made:
\begin{cor}\label{corollary-1}
If $S(N+1)\geq 0$, $R(k)>0$, and $\varPsi_k\geq 0$ hold for any $k\in\mathcal{Q}$,
then the LQ-optimal control problem \eqref{3} admits a unique optimal control and that control is in a state feedback form
$u^*(k)=-\mathcal{R}(k)^{-1}\mathcal{G}(k) x(k)$, $k\in\mathcal{Q}$,
and the optimal cost value $\mathcal{J}(x_0)=\langle P(0) x_0,x_0\rangle_H$.
\end{cor}
\begin{exmp}\label{example-1}
Convolution operators have widespread applications in signal processing, primarily used for signal filtering.
When it is applied to a signal, the Gaussian convolution operator convolves the signal with a Gaussian kernel
to reduce the high-frequency components of the signal and smooth its variations,
and thereby producing a new signal output (see, \citet{getreuer2013survey, spanias2007audio}).
For examples, \citet{ghalyan2018gaussian} employed a Gaussian smoothing filter to remove noise from the
electromyographic signals. \citet{kim2013applications} applied convolution operators into digital image signals with filters
to perform tasks such as edge extraction and reduction of unwanted noise.

Let $\mathbb{R}^1$ be one-dimensional real Euclidean space
and $L^2(-\infty,+\infty)=\{f|\Arrowvert f\Arrowvert_{L^2(-\infty,+\infty)}
=(\int_{-\infty}^{+\infty}|f(\tau)|^2\mathrm{d}\tau)^{\frac{1}{2}}\\<+\infty\}$.
By applying convolution operators to process and control audio signals,
a model for handling the interfering audio signals can be constructed as follows:
\begin{equation*}
	\begin{split}
		&X^{k+1}(t)=A(k)X^k(t)+B(k)u(k),\ k\in \mathcal{Q},\\
		&X^0(t)=X_0(t),
	\end{split}
\end{equation*}
where $X^k(t)\in L^2(-\infty,+\infty)$ represents the state of the $k$-th audio signal at time $t$.
Note that $t=0$ is not an absolute time point, but is defined relative to a specific event.
For instance, when processing an audio signal, a particular event,
such as the start of a sound or the generation of a sound wave, may be defined as $t=0$. To describe the state of the signal before this event, negative time values are used.
For any $k\in \mathcal{Q}$, $A(k)$ is a convolution operator with one-dimensional Gaussian kernel,
i.e., for any $f(t)\in  L^2(-\infty,+\infty)$,
\begin{equation*}
          A(k)f(t)=\int_{-\infty}^{+\infty}f(\tau)\frac{e^{-\frac{(t-\tau)^2}{2}}}{\sqrt{2\pi}}\mathrm{d}\tau,\ k\in \mathcal{Q}.
\end{equation*}
For any $k\in \mathcal{Q}$, $B(k)$ stands for an input operator that controls
the influence of the input $u(k)$ on the audio signal in the form of sound waves.
That is, for any $u(k)\in  \mathbb{R}^1$,
\begin{equation*}
	\begin{split}
            B(k)u(k)=\left\{\begin{aligned}
			&u(k) \ cos(\kappa \gamma -\omega t+\phi),\ &t\in[-1, 1],\\
			&0,\ &otherwise,
		\end{aligned}
		\right.
	\end{split}
\end{equation*}
where $u(k)$ is the control input, which can usually be regarded as the amplitude intensity of the signal wave,
and $\kappa$, $\gamma$, $\omega$ and $\phi$ represent the wave number, distance, angular frequency,
and initial phase of the sound wave, respectively.

For the interfering audio signals, we aim to control its state towards zero while minimizing the use of control inputs.
To achieve this, the following quadratic cost functional is introduced:
\begin{equation*}
	\begin{split}
		J(X_0(t),u)=& E\bigg[\sum_{k=0}^{N}\Big(\langle M(k)X^k,X^k \rangle_{ L^2(-\infty,+\infty) }\\
		&\ \ \ \ \ \ \ \ \ +\langle R(k)u(k),u(k)\rangle_{\mathbb{R}^1}\Big)\bigg].
	\end{split}
\end{equation*}
In which, $\langle M(k)X^k,X^k \rangle_{ L^2(-\infty,+\infty) }$ refers to
the penalty for the departure between the system state and the target state,
and $\langle R(k)u^k,u^k\rangle_{\mathbb{R}^1}$ refers to the penalty for excessive control inputs.

Set the initial audio signal as a Gaussian-modulated audio signal,
i.e., $X_0(t)=e^{\frac{-t^2}{2}}$ (see Fig.\ref{fig1}(a)).
And let $N=1$, $\kappa=\pi $, $\gamma=2$, $\omega=0.1\pi$,
$\phi=0$, $M(k)\equiv10I_{L^2(-\infty,+\infty)}$ and $R(k)\equiv 1$ ($k=0,1$).
By Corollary \ref{corollary-1}, one can acquire the optimal cost value $J(X_0(t),u^*)=24.052$,
and the optimal control $u^*=\{u^*(0),u^*(1)\}$ with $u^*(0)=-0.63$, $u^*(1)=0$.
Under the optimal control, the final signal state is shown in Fig.\ref{fig1}(b).

From Fig.\ref{fig1}(a), the initial audio signal exhibits a short-duration impulses or rapid variation,
resembling a sharp and brief sound.
After the control and convolution processing, the audio signal retains a smooth, gentle waveform
while experiencing a significant reduction in audio intensity, see Fig.\ref{fig1}(b).
\begin{figure}[ht]
	\centering
	\subfigure[]{
		\includegraphics[height=4cm]{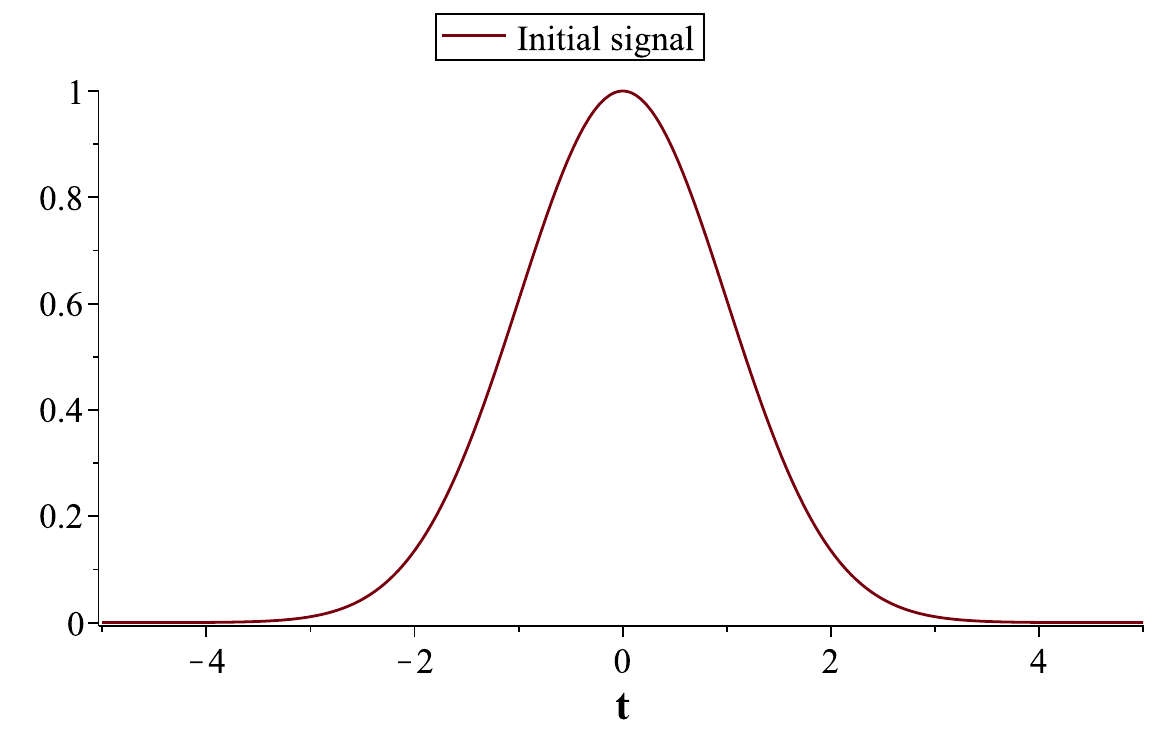}
	}
	\subfigure[]{
		\includegraphics[height=4cm]{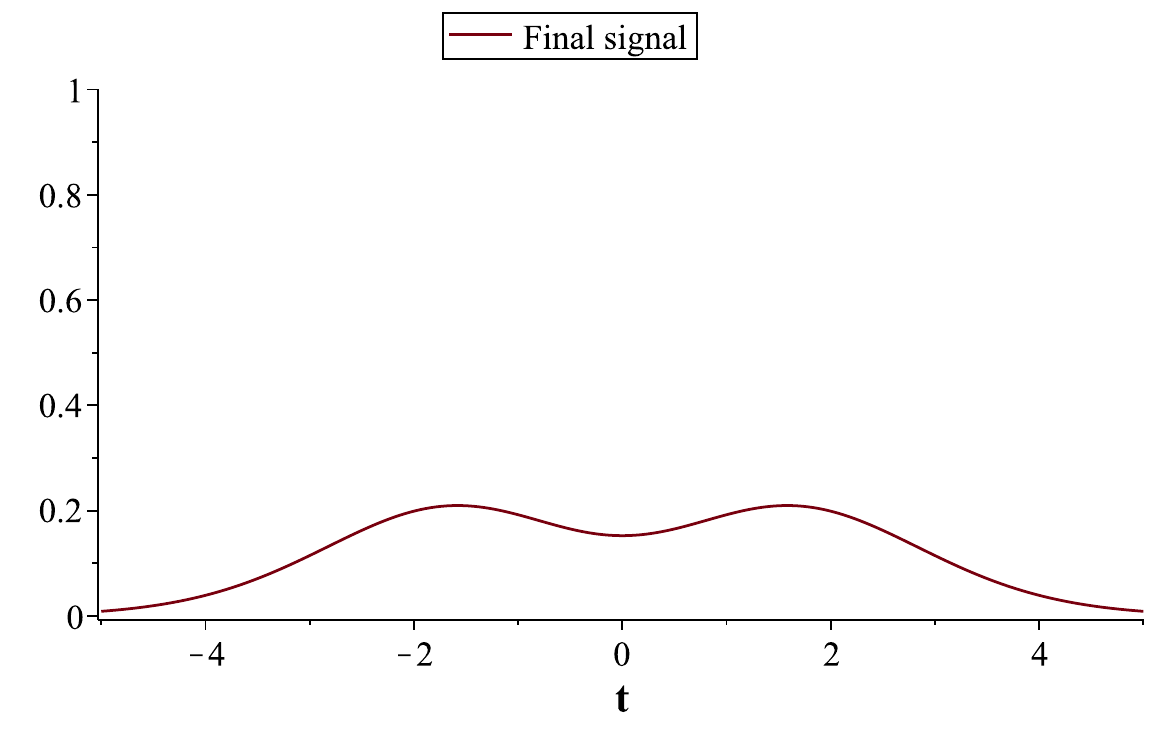}
	}
	\caption{$(a)$ Initial signal. $(b)$ Final signal.}	
    \label{fig1}
\end{figure}
\end{exmp}
\begin{rem}\label{remake-5}
In Example \ref{example-1}, for the optimal control $u^*=\{u^*(0),u^*(1)\}$,
note that $u^*(0)=-0.63$ does not imply that the audio signal's intensity is negative,
but rather represents a shift in the direction or phase of the waveform.
\end{rem}
\begin{exmp}\label{example-2}
In this example, one-dimensional temperature control system will be considered.
Denote $L^2[0,l]=\{f|\Arrowvert f\Arrowvert_{L^2[0,l]}
=(\int_{0}^{l}|f(\tau)|^2\mathrm{d}\tau)^{\frac{1}{2}}<+\infty\}$,
$H^1_0([0,l])=\{u\in L^2[0,l]| u^{'}\in L^2[0,l],\ u(0)=u(l)=0\}$,
and $H^2([0,l])=\{u\in L^2[0,l]| u^{'}, u^{''}\in L^2[0,l]\}$.
For all $T\in Dom(\Delta^{\alpha})=H^2([0,l])\cap H^1_0([0,l])$,
the Laplacian $\Delta^{\alpha}: Dom(\Delta^{\alpha})\mapsto L^2[0,l]$ is defined as
$(\Delta^{\alpha} T)(x)=\alpha\frac{ \partial^2T(x)}{\partial x^2}$.
Based on the theory of heat conduction and strongly continuous semigroups
(for examples, see \cite{2005Heat, Da1991, Flandoli1992}),
the following heat equation on $L^2[0,l]$:
\begin{equation*}
	\begin{split}
		\left\{\begin{aligned}
        &\frac{ \partial T(t, x)}{\partial t}=\alpha\frac{ \partial^2T(t, x)}{\partial x^2},\\
        &T(t, 0)=T(t, l)=0,\\
		&T(0, x)\in L^2[0,l],\ x\in [0,l]
		\end{aligned}
		\right.
	\end{split}
\end{equation*}
admits a solution $T(t, x)=S(t)T(0, x)\in L^2[0,l]$,
where the family of operators $S(t)=\{e^{t\Delta^{\alpha}}\}_{t\geq0}$ is a $C_0$-semigroup on $L^2[0,l]$
generated by the Laplacian $\Delta^{\alpha}$ with the zero Dirichlet boundary condition.
Now, by considering the representation of the solution and incorporating a control input $u$ at fixed time intervals $\tau$,
the following discrete-time linear control system is built:
\begin{equation*}
	\begin{split}
		&T^{k+1}(x)=A(k)T^k(x)+B(k)u(k),\ k\in \mathcal{Q},\\
		&T^0(x)=T_0(x)\in L^2[0,l],\ x\in [0,l],
	\end{split}
\end{equation*}
where $T^k(x)\in L^2[0,l]$ represents the temperature at time $k\tau$ in $x$,
and $u(k)\in \mathbb{R}^1$ is the control input vector (e.g., heating or cooling power) at time $k\tau$,
$\alpha$ is the coefficient of thermal conductivity.
For any $k\in \mathcal{Q}$, $A(k)=S(\tau)\in L(L^2[0,l])$ and $B(k)\in L(\mathbb{R}^1,L^2[0,l])$ respectively
describe the heat conduction characteristics and the impact of control inputs on the temperature.

Introduce a linear quadratic cost functional as follows:
\begin{equation*}
	\begin{split}
		J(T_0(x),u)=& E\bigg[\sum_{k=0}^{N}\Big(\langle M(k)T^k(x),T^k(x) \rangle_{L^2[0,l]}\\
		&+\langle R(k)u(k),u(k)\rangle_{\mathbb{R}^1}\Big)\\
        &+\langle S(N+1)T^{N+1}(x),T^{N+1}(x) \rangle_{L^2[0,l]}\bigg].
	\end{split}
\end{equation*}
In which, $M(k)$, $S(N+1)\in L(L^2[0,l])$ represent the penalties imposed on the deviation of the temperature state
from the target temperature $0^\circ C$ at time $k\tau$ and the terminal time $(N+1)\tau$, respectively.
$R(k)\in L(\mathbb{R}^1)$ denotes the penalty on the control input,
reflecting the energy consumption or implementation cost.

Set $N=2$, $l=1$, $\alpha=0.1 (W/(m\cdot K))$, $\tau=1 (s)$,
and the initial temperature $T_0(x)=60sin(\pi x)$,
as well as the control input $B(k)u(k)=u(k)x(l-x)$, $x\in[0,l]$.
Note that for any $k\in \mathcal{Q}$, the state vector $T^k(x)\in L^2[0,l]$ with the control input $B(k)u(k)$.
To better demonstrate the applicability and effectiveness of the proposed indefinite LQ-optimal control theory (Theorem \ref{theorem-1}),
we consider three representative cases with different cost function settings:

\noindent\textbf{Case 1:}
Let $M(k)\equiv 10I_{L^2[0,l]}$, $R(k)\equiv 1$ $(k=0,1,2)$, and $S(3)=10I_{L^2[0,l]}$.
This setting is common in practical applications,
and the goal is to achieve the desired state while avoiding excessive energy consumption.
In this case, the optimal control is $u^*=\{-33.3, -6.8, -2.5\}$, and the optimal cost value is $J(T_0(x),u^*)=22471$.

\noindent\textbf{Case 2:}
Let $M(k)\equiv 10I_{L^2[0,l]}$, $R(k)\equiv 0$ $(k=0,1,2)$, and $S(3)=10I_{L^2[0,l]}$.
It represents a scenario that energy consumption is not considered when achieving the target state.
In this case, the optimal control is $u^*=\{-122.4, -0.03, -0.01\}$, and the optimal cost value is $J(T_0(x),u^*)=18010$.

\noindent\textbf{Case 3:}
Let $M(k)\equiv 50I_{L^2[0,l]}$, $R(k)\equiv -1$ $(k=0,1,2)$, and $S(3)=50I_{L^2[0,l]}$.
It represents an extreme setting that excessive energy consumption is explicitly encouraged to achieve the desired state.
Such a scenario is rarely observed in practice.
In this setting, the optimal control is $u^*=\{-270.1, 6.2, 20.1\}$, and the optimal cost value is $J(T_0(x),u^*)=62243$.
\begin{figure}[ht]\label{fig2}
	\centering
	\subfigure[Case 1]{
		\includegraphics[height=4cm]{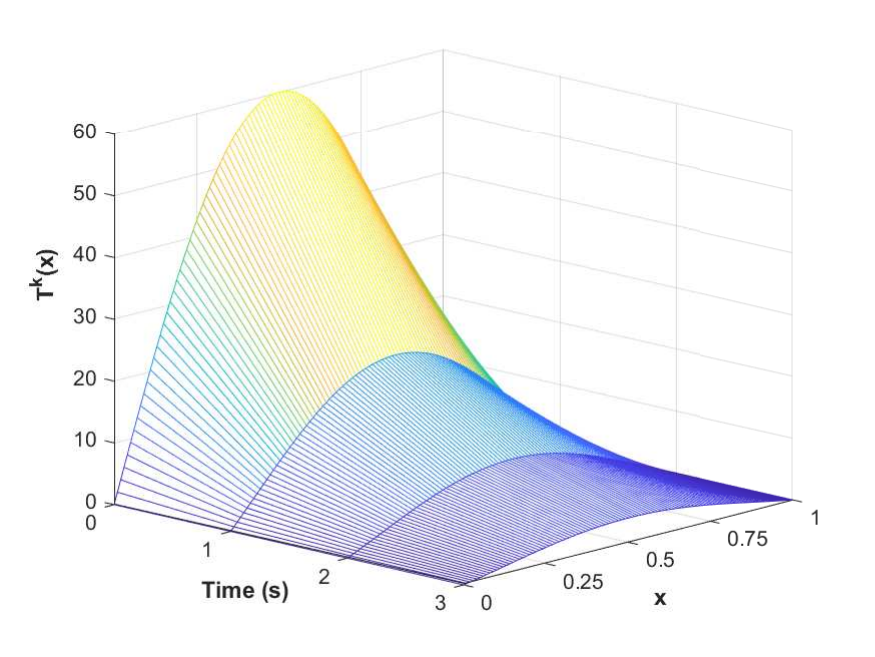}
	}
	\subfigure[Case 2]{
		\includegraphics[height=4cm]{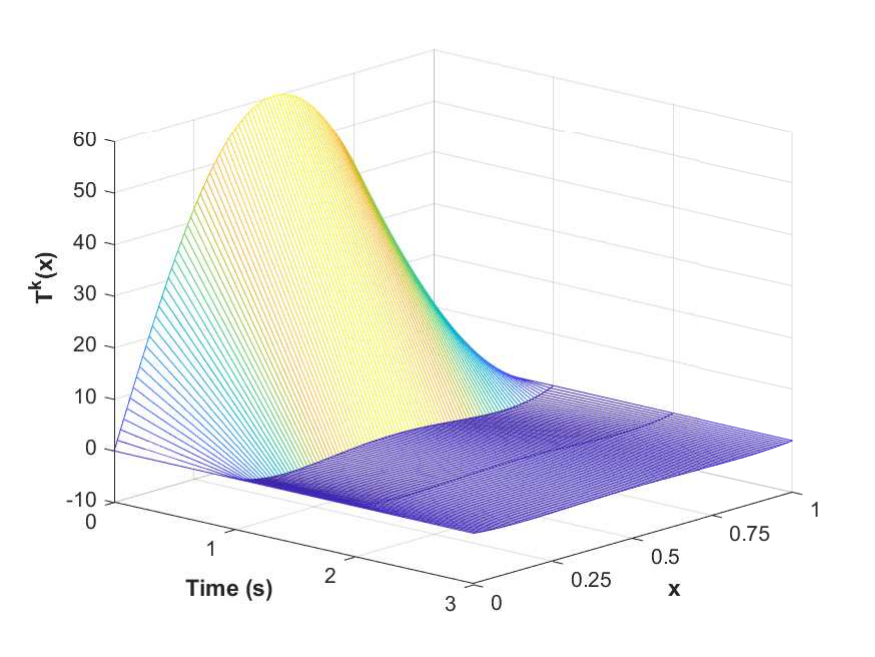}
	}
	\subfigure[Case 3]{
		\includegraphics[height=4cm]{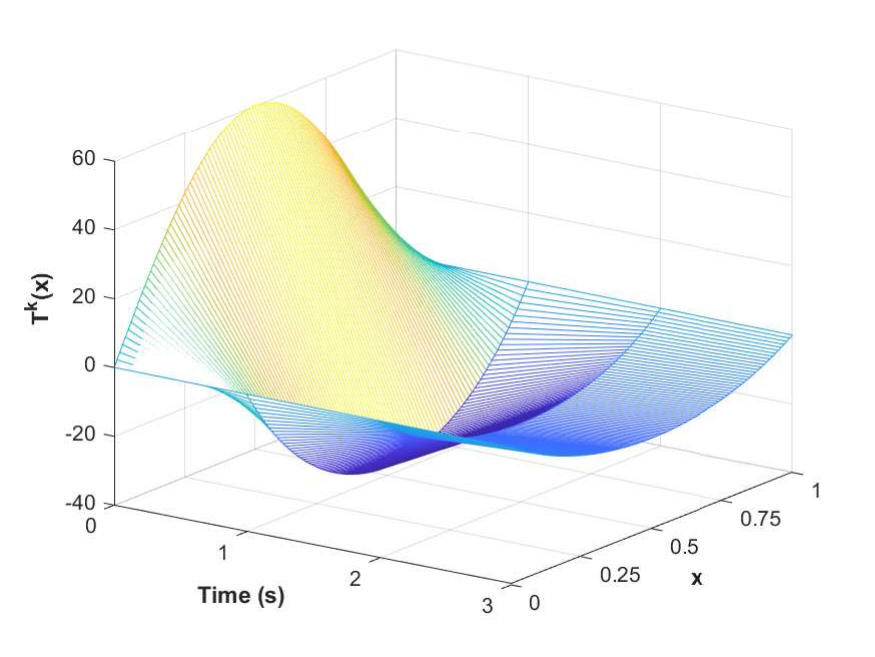}
	}
	\caption{The temperature at time $k$ in $x$.}	
\end{figure}

From Fig.\ref{fig2}(a), it can be observed that to avoid unnecessary energy loss,
under the optimal control, the temperature of the object does not plummet to $0^\circ C$,
but decreases gradually towards $0^\circ C$.
As shown in Fig.\ref{fig2}(b), since energy consumption is not considered,
the temperature rapidly drops to $0^\circ C$ under the optimal control.
Fig.\ref{fig2}(c) illustrates that although $R(k)$ takes a negative value,
due to the constraint imposed by the state cost term,
the optimal control can remain bounded and the temperature is prevented from deviating significantly from the target value $0^\circ C$.
\end{exmp}
\begin{rem}\label{remake-6}
In Example \ref{example-2}, based on the definition of the operator semigroup (see, \cite{Pazy1983}),
we can specify the form of the operator $A(k)$:
For $g(x)\in L^2[0,l]$,
	\begin{equation*}
        \begin{split}
       A(k)g(x)=\sum_{n=1}^{\infty}e^{-n^2\pi^2\alpha\tau}
                \langle g(x),\phi_n(x)\rangle_{L^2[0,l]}\phi_n(x),
	    \end{split}
    \end{equation*}
where $\phi_n(x)=\sqrt{\frac{2}{l}}sin(\frac{n\pi x}{l})$.
It can be verified that $A(k)$ is a self-adjoint and bounded linear operator.
Moreover, under the actions of the state operator $A(k)$ and the control operator $B(k)$,
the state $T^k(x)$ of the discrete-time linear control system satisfies the zero Dirichlet boundary condition at each time step.
\end{rem}
\begin{rem}	\label{remake-7}
It is straightforward to see that the available theoretical results
on the indefinite LQ-optimal control within the finite-dimensional setting (for example, see \cite{chen1998stochastic, 2016Chen})
are unable to handle the problems considered in Examples 1 and 2.
What's more, although substantial progress has been made in the studies of the LQ control problems under the infinite-dimensional framework,
most of the literature focuses on the continuous-time systems with convex cost functionals
(for example, see \cite{Flandoli1986, Da1988, 2022Hu, 2024Xu}).
However, the cost functionals arising in Case 2 and Case 3 of Example 2,
as well as the quadratic cost functional concerned with the upcoming $H_{\infty}$ analysis, are generally non-convex.
As a consequence, the existing results are insufficient to meet our needs,
and developing the indefinite LQ-optimal control theory for discrete-time infinite-dimensional systems is quite necessary,
which also serves as a foundation for the subsequent $H_{\infty}$ investigation.
\end{rem}

\section{Stochastic Bounded Real Lemma}

In this section, our goal is to analyze the disturbance attenuation for the following disturbed system:
\begin{equation}\label{8}
    \left\{\begin{aligned}
		&x(k+1)=A(k)x(k)+B_1(k)v(k)+(C(k)x(k)\\
		&\ \ \ \ \ \ \ \ \ \ \ \ \ \ \ \ +D_1(k)v(k))\omega(k),\\
		&z(k)=\overline{C}(k)x(k)+\overline{D}(k)v(k),\\
		&x(0)=x_0\in H,
	\end{aligned}
	\right.
\end{equation}
where $v(k)\in V$ and $z(k)\in Z$ ($k\in \mathcal{Q}$) stand for the exogenous disturbance signal and controlled output respectively,
$B_1(k)$, $D_1(k) \in L(V,H)$, $\overline{C}(k)\in L(H,Z)$ and $\overline{D}(k)\in L(V,Z)$,
and other symbols are defined as before.

Throughout this section, the following assumption is always satisfied:
\begin{assum}\label{assumption-1}
For any $k\in \mathcal{Q}$, the operators $\overline{C}(k)$ and $\overline{D}(k)$ are uncorrelated,
i.e., $\overline{D}(k)^*\overline{C}(k)=0$.
\end{assum}
\begin{defn}\label{definition-4}
For all $v=\{v(k),\ k\in \mathcal{Q}\}\in L^2_{V^{\mathcal{Q}}}$,
define the perturbation operator $\mathbb{L}:\ L^2_{V^{\mathcal{Q}}} \mapsto L^2_{Z^{\mathcal{Q}}}$ associated with \eqref{8}:
$(\mathbb{L}(v))(k)= \overline{C}(k)x(k)+\overline{D}(k)v(k),\ k\in \mathcal{Q},$
and its $H_{\infty}$ norm is
	\begin{equation*}
		\begin{split}
			\Arrowvert \mathbb{L} \Arrowvert&=\sup_{\substack{v\in L^2_{V^{\mathcal{Q}}},\\ v\neq 0,x_0=0}} \frac{\Arrowvert \mathbb{L}(v) \Arrowvert_{L^2_{Z^{\mathcal{Q}}}}}{\Arrowvert  v \Arrowvert_{L^2_{V^{\mathcal{Q}}}}}\\
			&=\sup_{\substack{v\in L^2_{V^{\mathcal{Q}}},\\ v\neq 0,x_0=0}} \frac{(\sum_{k=0}^{N} E[\Arrowvert \overline{C}(k)x(k)+\overline{D}(k)v(k)
\Arrowvert^2_{Z}])^{\frac{1}{2}}}{(\sum_{k=0}^N E[\Arrowvert v(k) \Arrowvert^2_{V}])^{\frac{1}{2}}}.
		\end{split}
	\end{equation*}
\end{defn}
The norm $\Arrowvert \mathbb{L} \Arrowvert$ quantifies the worst possible impact
that the disturbance $v$ might potentially inflict on the controlled output $z$ of \eqref{8}.
Introduce
    \begin{equation*}
		\begin{split}
            &\Lambda(x(k), v(k), \overline{C}(k), \overline{D}(k))\\
            &=\langle \overline{C}(k)^*\overline{C}(k)x(k),x(k)\rangle_{H}+\langle \overline{D}(k)^*\overline{D}(k)v(k),v(k)\rangle_{V}.
        \end{split}
	\end{equation*}
Then, on the basis of the properties of the inner product and Assumption \ref{assumption-1},
$\Arrowvert \mathbb{L} \Arrowvert$ can be equivalently expressed by
\begin{equation*}
		\Arrowvert \mathbb{L} \Arrowvert=
		 \sup_{\substack{v\in L^2_{V^{\mathcal{Q}}},\\v\neq 0,x_0=0}} \frac{(\sum_{k=0}^{N} E[\Lambda(x(k), v(k), \overline{C}(k), \overline{D}(k))])^{\frac{1}{2}}}{(\sum_{k=0}^N E[\langle v(k),v(k)\rangle_{V}])^{\frac{1}{2}}}.
\end{equation*}
Here, $x_0=0$ means that $x_0$ is set to be the zero element on Hilbert space $H$.
$v=\{v(k),\ k\in \mathcal{Q}\}\in L^2_{V^{\mathcal{Q}}}$ with $v\neq 0$ means that
there exists $k_0\in \mathcal{Q}$, $v(k_0)\in V$ and $v(k_0)\neq 0$.

Given $\gamma>0$, define the quadratic cost functional over a finite horizon as follows:
\begin{equation*}
	\begin{split}
		J_{\gamma}(0,v)=& E\bigg[\sum_{k=0}^{N}\langle -\overline{C}(k)^*\overline{C}(k)x(k),x(k) \rangle_H\\
		&\ \ \ \ +\langle (\gamma^2 I_V-\overline{D}(k)^*\overline{D}(k)) v(k),v(k)\rangle_V\bigg],
	\end{split}
\end{equation*}
where $x(\cdot)$ is the solution to \eqref{8} with $x_0=0$ and $v\in L^2_{V^{\mathcal{Q}}}$.

The following proposition connects the $H_{\infty}$ norm $\Arrowvert \mathbb{L} \Arrowvert$
with the quadratic cost functional $J_{\gamma}(0,v)$, and its proof is simple and therefore omitted.
\begin{prop}\label{proposition-2}
Given $\gamma>0$. For system \eqref{8} with $x_0=0$, $\Arrowvert \mathbb{L} \Arrowvert < \gamma$ iff
$\mathcal{J}_{\gamma}(0)=\inf_{v \in L^2_{V^{\mathcal{Q}}}, v\neq 0}J_{\gamma}(0,v)>0.$
\end{prop}

Furthermore, based on Proposition \ref{proposition-2}, some properties about the operator family
$\{\gamma^2 I_V-\overline{D}(k)^*\overline{D}(k), \\ k\in\mathcal{Q}\}$ can be gained.
\begin{prop}\label{proposition-3}
Given $\gamma>0$. For system \eqref{8} with $x_0=0$, if $\Arrowvert \mathbb{L} \Arrowvert<\gamma$ is satisfied,
then the operator family $\{\gamma^2 I_V-\overline{D}(k)^*\overline{D}(k), \ k\in\mathcal{Q}\}$ is uniformly positive.
\end{prop}
\begin{pf}
If the conclusion is not true, then there exists $k_0 \in \mathcal{Q} $ and $\widetilde{v}(k_0)\in V$, $\widetilde{v}(k_0)\neq0$
such that $\langle (\gamma^2 I_V-\overline{D}(k_0)^*\overline{D}(k_0))\widetilde{v}(k_0),\widetilde{v}(k_0)\rangle_V \leq 0$.
Take a disturbance signal $v^\Delta=\{v^\Delta(k),\ k\in \mathcal{Q}\}\in L^2_{V^{\mathcal{Q}}}$:
	\begin{equation*}
		v^\Delta(k)=\left\{\begin{aligned}
			&0, & k\neq k_0, \\
			&\widetilde{v}(k_0),  &k=k_0.
		\end{aligned}\right.
	\end{equation*}
Thus,
$J_{\gamma}(0,v^\Delta)=E[\sum_{k=k_0}^{N}\langle -\overline{C}(k)^*\overline{C}(k)x(k),x(k) \rangle_H
+\langle (\gamma^2 I_V-\overline{D}(k_0)^*\overline{D}(k_0)) \widetilde{v}(k_0),\widetilde{v}(k_0)\rangle_V]\leq0.$
So, we have that $\mathcal{J}_{\gamma}(0)\leq 0$.
By Proposition \ref{proposition-2}, it infers that $\Arrowvert \mathbb{L} \Arrowvert \geq \gamma$,
which contradicts with $\Arrowvert \mathbb{L} \Arrowvert<\gamma$.
The proof of Proposition \ref{proposition-3} is ended.
\end{pf}

Next, set $F=\{ F(k),\ k\in\mathcal{Q}\}\in L^2_{L(H,U)^\mathcal{Q}}$.
Given $\gamma>0$, define the operator $T^{\gamma}_k:\varepsilon(H) \mapsto \varepsilon(H)\ (k\in\mathcal{Q})$ as follows:
\begin{equation}\label{11}
	\begin{aligned}
		T^{\gamma}_k(X)=&A(k)^*XA(k)+C(k)^*XC(k)-\overline{C}(k)^*\overline{C}(k)\\
		&+[A(k)^*XB_1(k)+C(k)^*XD_1(k)]F(k)\\
		&+F(k)^*[B_1(k)^*XA(k)+D_1(k)^*XC(k)]\\
		&+F(k)^*[\gamma^2 I_V-\overline{D}(k)^*\overline{D}(k)+B_1(k)^*XB_1(k)\\
		&+D_1^*(k) XD_1(k)]F(k).
	\end{aligned}
\end{equation}
Consider the following backward operator equation:
\begin{equation}\label{12}
	\left\{\begin{aligned}
		&Y^{\gamma}(k)=T^{\gamma}_k(Y^{\gamma}(k+1)),\ k\in\mathcal{Q},\\
		&Y^{\gamma}(N+1)=0.
	\end{aligned}
	\right.
\end{equation}
Clearly, the solution to equation \eqref{12},
denoted as $Y^{\gamma}=\{Y^{\gamma}(k),\ k\in\mathcal{Q}^{+}\}\in L^2_{\varepsilon(H)^{\mathcal{Q}^{+}}}$,
can be acquired by the backward iteration.

In order to apply the backward operator equation \eqref{12} to handle the $H_{\infty}$ analysis of system \eqref{8},
we introduce the following matrix operators:
For any $X_1\in H$, $X_2\in V$ and $\mathcal{A}\in L(H)$, $\mathcal{B}\in L(H,V)$, define
\begin{equation*}
	\begin{split}
		&\left( \begin{matrix}
			\mathcal{A}\\ \mathcal{B}
		\end{matrix} \right) :H \mapsto H \times V, \ \
		\left( \begin{matrix}
			\mathcal{A}\\ \mathcal{B}
		\end{matrix} \right)^* :H\times V \mapsto H , \\
		&\left( \begin{matrix}
			\mathcal{A}\\ \mathcal{B}
		\end{matrix} \right)(X_1)= (\mathcal{A}X_1,\mathcal{B}X_1), \\
		&\left( \begin{matrix}
			\mathcal{A}\\ \mathcal{B}
		\end{matrix} \right)^*(X_1, X_2)=\mathcal{A}^*X_1+\mathcal{B}^*X_2.
	\end{split}
\end{equation*}
It is easy to see that the matrix operator
$\left( \begin{matrix}
	\mathcal{A}\\ \mathcal{B}
\end{matrix} \right)$ is also a bounded linear operator.
With the help of the matrix operators which have been well defined, the backward operator equation \eqref{12} can be rewritten as
\begin{equation}\label{14}
	\begin{split}
		\left( \begin{matrix}
			I_H\\ F(k)
		\end{matrix} \right)^*W(Y^{\gamma}(k),k)\left( \begin{matrix}
			I_H\\ F(k)
		\end{matrix} \right)=0,\ k\in\mathcal{Q},
	\end{split}
\end{equation}
where
\begin{equation*}
	\begin{split}
		&W(Y^{\gamma}(k),k)=\\
		&\ \  \ \left( \begin{matrix}
			\pi_1(Y^{\gamma}(k+1),k)-Y^{\gamma}(k) & \pi_2(Y^{\gamma}(k+1),k)^{*}\\
			\pi_2(Y^{\gamma}(k+1),k) & \pi^{\gamma}_3(Y^{\gamma}(k+1),k)
		\end{matrix} \right)
    \end{split}
\end{equation*}
with
\begin{equation*}
	\begin{split}
		&\pi_1(X,k)=A(k)^*XA(k)+C(k)^*XC(k)-\overline{C}(k)^*\overline{C}(k),\\
		&\pi_2(X,k)=B_1(k)^*XA(k)+D_1(k)^*XC(k),\\
		&\pi^{\gamma}_3(X,k)=\gamma^2 I_V-\overline{D}(k)^*\overline{D}(k)+B_1(k)^*XB_1(k)\\
		&\ \ \ \  \ \ \ \ \ \ \ \ \ \ \ \ +D_1(k)^*XD_1(k).
	\end{split}
\end{equation*}
\begin{lem}\label{lemma-3}
Fix $\gamma>0$ and let $F=\{ F(k),\ k\in\mathcal{Q}\}\in L^2_{L(H,U)^\mathcal{Q}}$.
If system \eqref{8} with $x_0=0$ satisfies $\Arrowvert \mathbb{L} \Arrowvert<\gamma$,
then the backward operator equation \eqref{12} admits a solution $\{Y^{\gamma}(k), \ k\in \mathcal{Q}^{+}\}$
such that $\pi^{\gamma}_3(Y^{\gamma}(k+1),k)>0$ holds for any $k\in \mathcal{Q}$.
\end{lem}
\begin{pf}
Firstly, we will prove that under the assumption, for any $k\in \mathcal{Q}$, $\pi^{\gamma}_3(Y^{\gamma}(k+1),k)\geq0$ can be deduced.
If it is false, then there exists $k_0\in \mathcal{Q}$ and $\widetilde{v}(k_0) \in V$, $\widetilde{v}(k_0)\neq 0$
such that $ \langle \pi^{\gamma}_3(Y^{\gamma}(k_0+1),k_0)\widetilde{v}(k_0),\widetilde{v}(k_0)\rangle_V<0$.
For $k\in \mathcal{Q}$, set
    \begin{equation}\label{1pq5}
         v^\Lambda(k)=F(k)x(k)+\hat{v}(k),
    \end{equation}
where
	\begin{equation*}
		 \hat{v}(k)=\left\{\begin{aligned}
			&0, & k\neq k_0, \\
			&\widetilde{v}(k_0),  &k=k_0.
		\end{aligned}\right.
	\end{equation*}
Obviously, $v^\Lambda=\{v^\Lambda(k),\ k\in \mathcal{Q}\}\in L^2_{V^{\mathcal{Q}}}$.
Additionally, if substituting $\hat{v}(k)$ into system \eqref{8} with $x_0=0$,
one can affirm that $x(k)=0$ for all $k\leq k_0$.
Consequently, from \eqref{1pq5}, it yields that $v^\Lambda(k_0)=F(k_0)x(k_0)+\hat{v}(k_0)=\widetilde{v}(k_0)\neq 0$,
which guarantees that for $k\in \mathcal{Q}$, $v^\Lambda(k)$ is not identic zero,
that is, $v^\Lambda\neq 0$.

In view of the fact that $\{Y^{\gamma}(k), \ k\in \mathcal{Q}\}$ is the solution to the backward operator equation,
it follows that
	\begin{equation*}
		\begin{split}
            &E\bigg[\sum_{k=0}^{N} (\langle W(Y^{\gamma}(k),k) (x(k),v^\Lambda(k)),(x(k),v^\Lambda(k))\rangle_{H \times V}\bigg]\\
			&\ \ \ -\langle -\overline{C}(k)^*\overline{C}(k)x(k),x(k) \rangle_H\\
			&\ \ \ -\langle (\gamma^2 I_V-\overline{D}(k)^*\overline{D}(k)) v^\Lambda(k),v^\Lambda(k)\rangle_V)\bigg]=0.
		\end{split}
	\end{equation*}
Based on this and \eqref{14}, we have that
	\begin{equation}\label{18}
		\begin{split}
			J_{\gamma}&(0,v)\\
			=&E\bigg[\sum_{k=0}^{N}\bigg(\bigg\langle \left( \begin{matrix}
				I_H\\ F(k)
			\end{matrix} \right)^*W(Y^{\gamma}(k),k)\left( \begin{matrix}
				I_H\\ F(k)
			\end{matrix} \right)x(k),x(k) \bigg\rangle_H\\
		    &\ \ \ +2\langle (\pi^{\gamma}_3(Y^{\gamma}(k+1),k)F(k)+\pi_2(Y^{\gamma}(k+1),k))\\
			&\ \ \ \cdot x(k),\hat{v}(k) \rangle_V+\langle \pi^{\gamma}_3(Y^{\gamma}(k+1),k)\hat{v}(k),\hat{v}(k)\rangle_V\bigg)\bigg].
		\end{split}
	\end{equation}
Since $\hat{v}(k)= 0$ for $k\neq k_0$ and $x(k)$ is the state response of system \eqref{8} with $x_0=0$,
\eqref{18} can be reduced to
$J_{\gamma}(0,v)=E[\langle \pi^{\gamma}_3(Y^{\gamma}(k_0+1),k_0)\widetilde{v}(k_0),\widetilde{v}(k_0)\rangle_V]<0,$
which contradicts with $\Arrowvert \mathbb{L} \Arrowvert<\gamma$,
$\pi^{\gamma}_3(Y^{\gamma}(k+1),k) \geq 0$ for any $k\in \mathcal{Q}$ is therefore concluded.
	
Secondly, we proceed to show that $\pi^{\gamma}_3(Y^{\gamma}(k+1),k)>0$ holds for any $k\in \mathcal{Q}$.
Set $0<\varepsilon_0<\gamma^2-\Arrowvert \mathbb{L} \Arrowvert^2$, $\tilde{\gamma}=(\gamma^2-\varepsilon_0)^{\frac{1}{2}}$.
Because of $\Arrowvert \mathbb{L} \Arrowvert<\tilde{\gamma}$,
following the discussion as above, one can deduce that $\pi^{\tilde{\gamma}}_3(Y^{\tilde{\gamma}}(k+1),k) \geq 0$ for any $k\in \mathcal{Q}$,
where $Y^{\tilde{\gamma}}=\{Y^{\tilde{\gamma}}(k),\ k\in \mathcal{Q}^{+}\}$ is a solution to the following backward operator equation:
	\begin{equation}\label{2po}
         \left\{\begin{aligned}
		  &Y^{\tilde{\gamma}}(k)=T^{\tilde{\gamma}}_k(Y^{\tilde{\gamma}}(k+1)),\ k\in\mathcal{Q},\\
		  &Y^{\tilde{\gamma}}(N+1)=0.
	     \end{aligned}\right.
    \end{equation}
Taking \eqref{11}, \eqref{12} and \eqref{2po} into account, it is valid that
	\begin{equation*}
		\begin{split}
			\left\{\begin{aligned}
				&Y^{\gamma}(k)-Y^{\tilde{\gamma}}(k)=[A(k)+B_1(k)F(k)]^*\\
				&\ \ \cdot[Y^{\gamma}(k+1)-Y^{\tilde{\gamma}}(k+1)][A(k)+B_1(k)F(k)]\\
				&\ \ +[C(k)+D_1(k)F(k)]^*[Y^{\gamma}(k+1)-Y^{\tilde{\gamma}}(k+1)]\\
				&\ \ \cdot[C(k)+D_1(k)F(k)]+\varepsilon_0F(k)^*F(k), \ k\in\mathcal{Q},\\
				&Y^{\gamma}(N+1)-Y^{\tilde{\gamma}}(N+1)=0.
			\end{aligned}
			\right.
		\end{split}
	\end{equation*}
We calculate by the backward iteration procedure and yield that
$Y^{\gamma}(k)-Y^{\tilde{\gamma}}(k)\geq 0,\ \forall k\in\mathcal{Q}.$
So, for any $k\in \mathcal{Q}$, it can be derived that
	\begin{equation*}
		\begin{split}
			&\pi^{\gamma}_3(Y^{\gamma}(k+1),k)\\
            &=B_1(k)^*Y^{\gamma}(k+1)B_1(k)+D_1(k)^*Y^{\gamma}(k+1)D_1(k)\\
			&\ \ \ \ +\gamma^2 I_V-\overline{D}(k)^*\overline{D}(k)\\
			&\geq B_1(k)^*Y^{\tilde{\gamma}}(k+1)B_1(k)+D_1(k)^*Y^{\tilde{\gamma}}(k+1)D_1(k)\\
			&\ \ \ \ +\tilde{\gamma}^2 I_V-\overline{D}(k)^*\overline{D}(k)+\varepsilon_0I_V\\
			&=\pi^{\tilde{\gamma}}_3(Y^{\tilde{\gamma}}(k+1))+\varepsilon_0I_V\\
			&>0.
		\end{split}
	\end{equation*}
Lemma \ref{lemma-3} is thus proved.
\end{pf}

Observing the proof of Lemma \ref{lemma-3}, it can be found that
if for $k\in \mathcal{Q}$, we set $v^\Lambda(k)=\hat{v}(k)$
(i.e., $F(k)=0$ in \eqref{1pq5}),
then Proposition \ref{proposition-3} is drawn.
Hence, Proposition \ref{proposition-3} may also be seen as a corollary of Lemma \ref{lemma-3}.

The next conclusion, called the finite-horizon stochastic bounded real lemma,
provides a sufficient and necessary condition for system \eqref{8}
to ensure its $H_{\infty}$ norm $\Arrowvert \mathbb{L} \Arrowvert$ below a specified level $\gamma >0$.
\begin{thm}\label{theorem-3}(Finite-horizon Stochastic Bounded Real Lemma)
Given $\gamma>0$, system \eqref{8} with $x_0=0$ satisfies $\Arrowvert \mathbb{L} \Arrowvert<\gamma$
iff the following backward operator equation:
	\begin{equation}\label{23}
		\left\{\begin{aligned}
			&Y^{\gamma}(k)=\pi_1(Y^{\gamma}(k+1),k)-\pi_2(Y^{\gamma}(k+1),k)^{*}\\
			&\ \ \cdot\pi^{\gamma}_3(Y^{\gamma}(k+1),k)^{-1}\pi_2(Y^{\gamma}(k+1),k),\ \forall k\in\mathcal{Q},\\
			&Y^{\gamma}(N+1)=0
		\end{aligned}
		\right.
	\end{equation}
admits a solution $\{ Y^{\gamma}(k),\ k\in\mathcal{Q}^{+}\}\in L^2_{\varepsilon(H)^{{\mathcal{Q}^{+}}}}$ such that
$\pi^{\gamma}_3(Y^{\gamma}(k+1),k)>0,\ \forall k\in \mathcal{Q}.$
\end{thm}
\begin{pf}
Sufficiency.
In \eqref{2}, take $M(k)=-\overline{C}(k)^*\overline{C}(k)$, $L(k)=0$, $R(k)=\gamma^2 I_V-\overline{D}(k)^*\overline{D}(k)$ and $S(N+1)=0$.
Applying Theorem \ref{theorem-1} to system \eqref{8}, we can directly obtain that $\mathcal{J}_{\gamma}(0)=0$,
and $v^*(k)=-\pi^{\gamma}_3(Y^{\gamma}(k+1),k)^{-1}\pi_2(Y^{\gamma}(k+1),k)x(k),\ k\in \mathcal{Q}$.
Now substituting $v^*(k)$ ($k\in \mathcal{Q}$) into system \eqref{8},
one has that $x(k)=0$ for any $k\in \mathcal{Q}$.
Therefore, $J_{\gamma}(0,v) >0$ whenever $v\neq 0$.
So, $\Arrowvert \mathbb{L} \Arrowvert<\gamma$ is concluded by Proposition \ref{proposition-2}.
	
Necessity.
Firstly, set $F=\{F(0)=0,F(1)=0,\cdots,F(N)=0\}$.
By Lemma \ref{lemma-3}, we have that $\pi^{\gamma}_3(Y^{\gamma}(N+1),N)>0$,
which implies that $-\pi^{\gamma}_3(Y^{\gamma}(N+1),N)^{-1}\pi_2(Y^{\gamma}(N+1),N)$ is well defined.
Then, substituting $\{F(0)=0,\cdots,F(N-1)=0,F(N)=-\pi^{\gamma}_3(Y^{\gamma}(N+1),N)^{-1}\pi_2(Y^{\gamma}(N+1),N)\}$ into \eqref{12},
the resulting equation is consistent with equation \eqref{23}.
By Lemma \ref{lemma-3} again, $Y^{\gamma}(N)$ exists and $\pi^{\gamma}_3(Y^{\gamma}(N),N-1)>0$.
Next, by constructing $\{F(0)=0,\cdots,F(N-2)=0,F(N-1)=-\pi^{\gamma}_3(Y^{\gamma}(N),N-1)^{-1}\pi_2(Y^{\gamma}(N),N-1), F(N)=-\pi^{\gamma}_3(Y^{\gamma}(N+1),N)^{-1}\pi_2(Y^{\gamma}(N+1),N)\}$
and using Lemma \ref{lemma-3},
$Y^{\gamma}(N-1)$ exists and $\pi^{\gamma}_3(Y^{\gamma}(N-1),N-2)>0$.
Continuing this process until $k=0$,
$Y^{\gamma}(k)$, $k=N-2,N-3,\cdots,0$ exist and $\pi^{\gamma}_3(Y^{\gamma}(k+1),k)>0$, $k=N-3,N-4,\cdots,0$.
Thereby, the proof of necessity is complete.
\end{pf}
\begin{rem}\label{remake-8}
Since $Y^{\gamma}(N)=-\overline{C}(N)^*\overline{C}(N)\leq 0$,
which suggests that the solution $\{ Y^{\gamma}(k),\ k\in\mathcal{Q}^{+}\}$ to the backward operator equation \eqref{23}
fails to be uniformly nonnegative,
the Schur's Complement Lemma cannot be applied when dealing with \eqref{23}.	
\end{rem}
\begin{rem}\label{remake-9}
It is worth noting that the extension of stochastic bounded real lemma
in the infinite-dimensional context is fundamentally nontrivial.
The $H_{\infty}$ analysis performed in this paper explicitly incorporates the structural features of bounded linear operators
and the stochastic evolution of the disturbed system.
In addition, the indefinite nature of the quadratic cost functional invalidates the convex cost functional-based LQ arguments
and necessitates a careful characterization of positivity and solvability conditions on the backward Riccati operator equations.
\end{rem}
\begin{exmp}\label{example-3}
Consider the discrete-time linear system with multiplicative noise and disturbance \eqref{8}.
Set $N=5$, and the spaces $H=Z=\ell^2$ and $V=\mathbb{R}^4$.
For any $k\in \mathcal{Q}$, let $\overline{C}(k)$ be the right shift operator,
and $\overline{D}(k)$ be the first unit filling operator for four-dimensional elements,
i.e., $\overline{D}(k)(a_1, a_2, a_3, a_4)=(a_1, 0, 0, 0, \cdots)$.
In addition, the operators $A(k)$, $B_1(k)$, $C(k)$ and $D_1(k)$\ ($k\in \mathcal{Q}$) are taken as follows:
\begin{equation*}
	\begin{split}
		k\ is\ &an\ odd\ number:\\
		&A(k)=C(k)\ are\ compression\ right\ shift\ operators:\\
		&\ \ \ A(k)(a_1, a_2, a_3,\cdots)=\frac{\sqrt{2}}{4}(0,a_1, a_2, a_3,\cdots),
    \end{split}
\end{equation*}
\begin{equation*}
	\begin{split}
		&B_1(k)=D_1(k)\ are\ compression\ filling\ operators:\\
		&\ \ \ B_1(k)(a_1, a_2, a_3, a_4)=\frac{\sqrt{2}}{4}(a_1, a_2, a_3, a_4, 0, \cdots);\\
		k\ is\ &an\ even\ number:\\
		&A(k)=C(k)\ are\ compression\ operators:\\
		&\ \ \ A(k)x=\frac{\sqrt{2}}{4}I_H x, \ x\in H,\\
		&B_1(k)=D_1(k)=0.
	\end{split}
\end{equation*}
Under the above settings, let us consider the backward operator equation \eqref{23}.
Notice that for any $k\in \mathcal{Q}$, $\pi^{\gamma}_3(Y^{\gamma}(k+1),k)\in L(V)$,
and the minimum point spectra $\rho_{min}$ of $\pi^{\gamma}_3(Y^{\gamma}(k+1),k)$ are
\begin{equation*}
	\begin{split}
		&\rho_{min}(\pi^{\gamma}_3(Y^{\gamma}(k+1),k))=\gamma^2-1,\ k=0,2,4,5,\\
		&\rho_{min}(\pi^{\gamma}_3(Y^{\gamma}(4),3))=\gamma^2-\frac{9}{4},\\
		&\rho_{min}(\pi^{\gamma}_3(Y^{\gamma}(2),1))=\gamma^2-\frac{41}{16}-\frac{25}{64(\gamma^2-\frac{5}{4})},\\
	\end{split}
\end{equation*}
Therefore, ensuring $\{\pi^{\gamma}_3(Y^{\gamma}(k+1),k),\ k\in \mathcal{Q}\}$ to be uniformly positive is equivalent to
$\rho_{min}(\pi^{\gamma}_3(Y^{\gamma}(k+1),k))>0,\ k=0,1,\cdots,5,$
which requires $\gamma^2>\frac{45}{16}$.
Thus, according to the finite-horizon stochastic bounded real lemma (Theorem \ref{theorem-3}),
system \eqref{8} with initial value $x_0=0$ satisfies $\Arrowvert \mathbb{L} \Arrowvert=\frac{3\sqrt{5}}{4}$.
\end{exmp}

\section{Game-based $H_{\infty}$ and $H_2/H_{\infty}$ Control Designs}

In this section, the Nash equilibrium problem is studied firstly,
which evolves a unified approach for the $H_{\infty}$ and $H_2/H_{\infty}$ control designs.
Let us consider the following discrete-time stochastic system with two inputs:
\begin{equation}\label{24}
	\left\{\begin{aligned}
		&x(k+1)=A(k)x(k)+B_1(k)v(k)+B_2(k)u(k)\\
		&\ \ \ +(C(k)x(k)+D_1(k)v(k)+D_2(k)u(k))\omega(k),\\
		&z(k)=\overline{C}(k)x(k)+\overline{G}(k)u(k),\\
		&x(0)=x_0\in H,
	\end{aligned}
	\right.
\end{equation}
where $v(k)\in V$ and $u(k)\in U$ ($k\in \mathcal{Q}$) are the inputs manipulated by two different players,
$B_1(k)$, $D_1(k)\in L(V,H)$, $B_2(k)$, $D_2(k)\in L(U,H)$,
$\overline{G}(k)\in L(U,Z)$, and other symbols are defined as before.

This section adopts the following basic assumption:
\begin{assum}\label{assumption-2}
For any $k\in \mathcal{Q}$, the operators $\overline{C}(k)$ and $\overline{G}(k)$ satisfy
$\overline{G}(k)^*\overline{C}(k)=0$ and $\overline{G}(k)^*\overline{G}(k)=I_U$.
\end{assum}

To investigate the Nash equilibrium problem associated with system \eqref{24},
we first introduce two parameterized quadratic performance indices as follows:
\begin{eqnarray}
	J_1^N(x_0,u,v):= E\bigg[\sum_{k=0}^{N}\Big(\gamma^2\Arrowvert v(k) \Arrowvert_V^2-\Arrowvert z(k) \Arrowvert_Z^2\Big)\bigg],\label{25}\\
	J_2^N(x_0,u,v):= E\bigg[\sum_{k=0}^{N}\Big(\Arrowvert z(k) \Arrowvert_Z^2-\rho^2\Arrowvert v(k) \Arrowvert_V^2\Big)\bigg],\label{26}
\end{eqnarray}
where $u=\{u(k),\ k\in \mathcal{Q}\}\in L^2_{U^{\mathcal{Q}}}$, $v=\{v(k),\ k\in \mathcal{Q}\}\in L^2_{v^{\mathcal{Q}}}$,
and $0<\gamma<+\infty$ and $0\leq\rho<+\infty$ are real parameters.
\begin{defn}\label{definition-5}
The equilibrium strategy pair $(u^*,v^*)\in L^2_{U^{\mathcal{Q}}}\times L^2_{V^{\mathcal{Q}}}$ is called a Nash equilibrium,
if for any $(u,v)\in L^2_{U^{\mathcal{Q}}}\times L^2_{V^{\mathcal{Q}}}$,
	\begin{equation}\label{27}
		J_1^N(x_0,u^*,v^*)\leq J_1^N(x_0,u^*,v),
	\end{equation}
and 	
	\begin{equation}\label{28}
		J_2^N(x_0,u^*,v^*)\leq J_2^N(x_0,u,v^*)
	\end{equation}
are satisfied.
\end{defn}

Below, it is attempted to find a Nash equilibrium $(u^*,v^*)$ with the form of linear feedbacks.
Before proceeding, we present the cross-coupled operators as follows:
For $k\in \mathcal{Q}$ and $(X_1,X_2)\in L^2_{\varepsilon(H)^{{\mathcal{Q}^{+}}}}
\times L^2_{\varepsilon(H)^{{\mathcal{Q}^{+}}}}$,
\begin{equation*}
	\begin{split}
		&\mathcal{R}_1(k, X_1) =\gamma^2I_V+B_1(k)^*X_1(k+1)B_1(k)\\
        &\ \ \ \ \ \ \ \ \ \ \ \ \ \ \ \ \ \ \ +D_1(k)^*X_1(k+1)D_1(k),\\
        &\mathcal{R}_2(k,X_2) =I_U+B_2(k)^*X_2(k+1)B_2(k)\\
		&\ \ \ \ \ \ \ \ \ \ \ \ \ \ \ \ \ \ \ +D_2(k)^*X_2(k+1)D_2(k),\\
        &\mathcal{G}_1(k,X_1) = B_1(k)^*X_1(k+1)[A(k)+B_2(k)K_2(k,X_2)]\\
        &\ \ \ \ \ \ \ \ \ \ \ \ \ \ +D_1(k)^*X_1(k+1)[C(k)+D_2(k)K_2(k,X_2)],\\
        &\mathcal{G}_2(k,X_2) = B_2(k)^*X_2(k+1)[A(k)+B_1(k)K_1(k,X_1)]\\
		&\ \ \ \ \ \ \ \ \ \ \ \ \ \ +D_2(k)^*X_2(k+1)[C(k)+D_1(k)K_1(k,X_1)],\\
    	&K_1(k,X_1)=-\mathcal{R}_1(k,X_1)^{-1}\mathcal{G}_1(k,X_1),\\
		&K_2(k,X_2)=-\mathcal{R}_2(k,X_2)^{-1}\mathcal{G}_2(k,X_2).
	\end{split}
\end{equation*}
Based on the above operators, define the following coupled backward Riccati operator equations:
\begin{equation}\label{29}
	\left\{\begin{aligned}
		&P_1(k)=[A(k)+B_2(k)K_2(k,P_2)]^*P_1(k+1)\\
		&\ \ \ \ \ \ \ \ \ \ \ \ \ \cdot [A(k)+B_2(k)K_2(k,P_2)]\\
		&\ \ \ \ \ \ \ \ \ \ \ \ \ +[C(k)+D_2(k) K_2(k,P_2)]^*P_1(k+1)\\
        &\ \ \ \ \ \ \ \ \ \ \ \ \ \cdot [C(k)+D_2(k)K_2(k,P_2)]\\
		&\ \ \ \ \ \ \ \ \ \ \ \ \ -K_2(k,P_2)^*K_2(k,P_2)-\overline{C}(k)^*\overline{C}(k)\\
		&\ \ \ \ \ \ \ \ \ \ \ \ \ -\mathcal{G}_1(k,P_1)^*\mathcal{R}_1(k,P_1)^{-1}\mathcal{G}_1(k,P_1),\\
        &\mathcal{R}_1(k,P_1)^{-1}>0,\ k\in\mathcal{Q},\\
		&P_1(N+1)=0,
	\end{aligned}
	\right.
\end{equation}
\begin{equation}\label{30}
	\left\{\begin{aligned}
		&P_2(k)=[A(k)+B_1(k)K_1(k,P_1)]^*P_2(k+1)\\
		&\ \ \ \ \ \ \ \ \ \ \ \ \ \cdot [A(k)+B_1(k)K_1(k,P_1)]\\
		&\ \ \ \ \ \ \ \ \ \ \ \ \ +[C(k)+D_1(k)K_1(k,P_1)]^*P_2(k+1)\\
        &\ \ \ \ \ \ \ \ \ \ \ \ \ \cdot [C(k)+D_1(k)K_1(k,P_1)]\\
		&\ \ \ \ \ \ \ \ \ \ \ \ \ -\rho^2K_1(k,P_1)^*K_1(k,P_1)+\overline{C}(k)^*\overline{C}(k)\\
		&\ \ \ \ \ \ \ \ \ \ \ \ \ -\mathcal{G}_2(k,P_2)^*\mathcal{R}_2(k,P_2)^{-1}\mathcal{G}_2(k,P_2),\\
        &\mathcal{R}_2(k,P_2)^{-1}>0,\ k\in\mathcal{Q},\\
		&P_2(N+1)=0.
	\end{aligned}
	\right.
\end{equation}
\begin{thm}\label{theorem-4}
For system \eqref{24}, the two-player Nash game problem \eqref{27}-\eqref{28} is solved via a linear feedback Nash equilibrium
iff the coupled backward Riccati operator equations \eqref{29}-\eqref{30}
admit a pair of solutions $(P_1,P_2)\in L^2_{\varepsilon(H)^{{\mathcal{Q}^{+}}}}\times L^2_{\varepsilon(H)^{{\mathcal{Q}^{+}}}}$.
And the optimal Nash strategy pair $(u^*,v^*)$ is given by
$u^*=\{u^*(k)=K_2(k,P_2)x(k),\ k\in \mathcal{Q}\}$
and $v^*=\{v^*(k)=K_1(k,P_1)x(k),\ k\in \mathcal{Q}\},$
where $x(k)$ ($k\in \mathcal{Q}$) is the current state of
\begin{equation*}
	\left\{\begin{aligned}
		&x(k+1)=(A(k)+B_1(k)K_1(k,P_1)+B_2(k)\\
&\ \ \ \cdot K_2(k,P_2))x(k)+(C(k)+D_1(k)K_1(k,P_1)\\
&\ \ \ +D_2(k)K_2(k,P_2))x(k)\omega(k),\\
		&x(0)=x_0\in H.
	\end{aligned}
	\right.
\end{equation*}
In this case, the cost values incurred by $(u^*,v^*)$ are
$J_1^N(x_0,u^*,v^*)=\langle P_1(0)x_0,x_0 \rangle_H$
and $J_2^N(x_0,u^*,v^*)=\langle P_2(0)x_0,x_0 \rangle_H.$
\end{thm}
\begin{pf}
See Appendix.
\end{pf}

By selecting specific values for the parameters $\gamma$ and $\rho$ in \eqref{25}-\eqref{26},
Theorem \ref{theorem-4} can be applied to solve the $H_{\infty}$ and $H_2/ H_{\infty}$ control problems of system \eqref{24}.

\noindent\textbf{Case I: $0<\gamma=\rho<+\infty$}

In this case, the Nash game problem \eqref{25}-\eqref{26} reduces to a zero-sum game, that is,
\begin{equation}\label{123po}
	\begin{aligned}
		J_1^N(x_0,u,v)+ J_2^N(x_0,u,v)=0.
	\end{aligned}
\end{equation}
Assume that $P_1(k)$ and $P_2(k)$ ($k\in\mathcal{Q}^{+}$) are the solutions to
the coupled backward Riccati operator equations \eqref{29} and \eqref{30}, respectively.
By the completing the square technique, one can calculate that
\begin{equation*}
	\begin{split}
		&P_1(k)+P_2(k)\\
        &\ \ =[A(k)+B_1(k)K_1(k,P_1)+B_2(k)K_2(k,P_2)]^*\\
		&\ \ \ \ \ \ \cdot [P_1(k+1)+P_2(k+1)][A(k)+B_1(k)K_1(k,P_1)\\
		&\ \ \ \ \ \ +B_2(k)K_2(k,P_2)]+[C(k)+D_1(k)K_1(k,P_1)\\
		&\ \ \ \ \ \ +D_2(k)K_2(k,P_2)]^*[P_1(k+1)+P_2(k+1)]\\
		&\ \ \ \ \ \ \cdot [C(k)+D_1(k)K_1(k,P_1)+D_2(k)K_2(k,P_2)].
	\end{split}
\end{equation*}
In accordance with $P_1(N+1)+P_2(N+1)=0$, then for any $k\in\mathcal{Q}$, $P_1(k)+P_2(k)=0$ is immediately received from above.
Let $P(k)=-P_1(k)=P_2(k)$, $k\in\mathcal{Q}$.
Then, the coupled backward Riccati operator equations \eqref{29} and \eqref{30} give that
\begin{equation}\label{41}
	\left\{\begin{aligned}
		&P(k)=[A(k)+B_2(k)K_2(k,P)]^*P(k+1)\\
		&\ \ \ \ \ \ \ \ \ \ \ \cdot [A(k)+B_2(k)K_2(k,P)] \\
        &\ \ \ \ \ \ \ \ \ \ \ +[C(k)+D_2(k)K_2(k,P)]^*\\
		&\ \ \ \ \ \ \ \ \ \ \ \cdot P(k+1)[C(k)+D_2(k)K_2(k,P)]\\
        &\ \ \ \ \ \ \ \ \ \ \ +K_2(k,P)^*K_2(k,P)+\overline{C}(k)^*\overline{C}(k)\\
		&\ \ \ \ \ \ \ \ \ \ \ +\mathcal{G}_1(k,-P)^*\mathcal{R}_1(k,-P)^{-1}\mathcal{G}_1(k,-P),\\
		&\mathcal{R}_1(k,-P)=\gamma^2I_V-B_1(k)^*P(k+1)B_1(k)\\
		&\ \ \ \ \ \ \ \ \ -D_1(k)^*P(k+1)D_1(k)>0,\\
        &\mathcal{R}_2(k,P)=I_U+B_2(k)^*P(k+1)B_2(k)\\
		&\ \ \ \ \ \ \ \ \ +D_2(k)^*P(k+1)D_2(k)>0,\ \forall k\in\mathcal{Q},\\
		&P(N+1)=0.
	\end{aligned}
	\right.
\end{equation}
Meanwhile, the optimal Nash strategy pair is
$u^*(k)=K_2(k,P)x(k)$, $v^*(k)=K_1(k,-P)x(k)$, $k\in\mathcal{Q}$,
where $P=\{P(k),\ k\in\mathcal{Q}^{+}\}$ is the solution to the backward Riccati operator equation \eqref{41}.
Furthermore, in light of Definition \ref{definition-5} and \eqref{123po}, one can work out that
$J_1^N(x_0,u,v^*)\leq J_1^N(x_0,u^*,v^*)\leq J_1^N(x_0,u^*,v).$

Now, we turn our attention to the $H_{\infty}$ control problem associated with system \eqref{24}.
In \eqref{24}, set $x_0=0$, and $u(\cdot)$ and $v(\cdot)$ are regarded as
the control input and the exogenous disturbance signal, respectively.
On the basis of the above analysis, together with Theorems \ref{theorem-3} and \ref{theorem-4},
the following theorem can be derived directly.
\begin{thm}($H_{\infty}$ Control)\label{theorem-5}
Consider system \eqref{24} with $x_0=0$.
Given $\gamma>0$, there exists a linear feedback control sequence $\{u^*(k),\ k\in\mathcal{Q}\}$ such that
the $H_{\infty}$ norm $\Arrowvert \mathbb{L}_{u^*} \Arrowvert$ of the corresponding closed-loop system
satisfies $\Arrowvert \mathbb{L}_{u^*} \Arrowvert<\gamma$
iff the backward Riccati operator equation \eqref{41} admits a solution
$P\in L^2_{\varepsilon(H)^{\mathcal{Q}}}$.
And, $u^*=\{u^*(k)=K_2(k,P)x(k),\ k\in\mathcal{Q}\}$ is the $H_{\infty}$ control,
where $P=\{P(k),\ k\in\mathcal{Q}^{+}\}$ is the solution to \eqref{41}.
\end{thm}

\noindent\textbf{Case II: $0<\gamma<+\infty$, $\rho=0$}

For system \eqref{24}, when it is imposed $\rho=0$ on $J_2^N(x_0,u,v)$, we have that
\begin{equation}\label{290}
	J_2^N(x_0,u,v)= E\bigg[\sum_{k=0}^{N}\Arrowvert z(k) \Arrowvert_Z^2\bigg].
\end{equation}

The following result concerning with the $H_2/H_{\infty}$ controller design can be yielded from Theorem \ref{theorem-4} straightforwardly.

\begin{thm}($H_2/H_{\infty}$ Control)\label{theorem-6}
Consider system \eqref{24} with $x_0=0$.
Given $\gamma>0$, there exists a pair of linear feedback equilibrium $(u^*(\cdot),v^*(\cdot))$ such that the following holds:\\
$(i)$ Substituting $u=u^*$ in \eqref{24}, the corresponding closed-loop system satisfies $\Arrowvert \mathbb{L}_{u^*} \Arrowvert<\gamma$;\\
$(ii)$ When the worst-case disturbance $v^*$ (see its definition in \cite{Limebeer1994nash}) is implemented,
$u^*$ minimizes $J_2^N(x_0,u,v)$ which is defined by \eqref{290}, i.e.,
$$J_2^N(x_0,u^*,v^*)=\inf_{u \in \mathcal{U}}J_2^N(x_0,u,v^*)$$
iff the coupled backward Riccati operator equations \eqref{29} and \eqref{30} with $\rho=0$
admit a pair of solutions $(P_1,P_2)\in L^2_{\varepsilon(H)^{{\mathcal{Q}^{+}}}}\times L^2_{\varepsilon(H)^{{\mathcal{Q}^{+}}}}$.
In this case, the control law is $u^*=\{u^*(k)=K_2(k,P_2)x(k),\ k\in\mathcal{Q}\}$
and the worst-case disturbance is $v^*=\{v^*(k)=K_1(k,P_1)x(k),\ k\in\mathcal{Q}\}$,
where $P_1=\{P_1(k),\ k\in\mathcal{Q}^{+}\}$ and $P_2=\{P_2(k),\ k\in\mathcal{Q}^{+}\}$
are the solutions to \eqref{29} and \eqref{30} with $\rho=0$, respectively.
In addition, $J_2^N(x_0,u^*,v^*)=\langle P_2(0)x_0,x_0 \rangle_H$.
\end{thm}
\begin{rem}\label{remake-10}
It must be pointed out that the $H_2$ control design cannot be carried out
by setting $\gamma\rightarrow+\infty$ and $\rho=0$ in Theorem \ref{theorem-4}
because the well-posedness is not ensured in this setting.
This is different from the finite-dimensional situation discussed in \cite{2013A, zhang2017stochastic}.
\end{rem}
At the end of this subsection, we present an example concerned with the two-player Nash game problem.
\begin{exmp}\label{example-4}
In \eqref{24}, set $N=1$, $H=Z=\ell^2$, $U=V=\mathbb{R}^1$.
For any $k\in \mathcal{Q}$, let $A(k)=C(k)=\frac{1}{2}I_H$, $B_1(k)=B_2(k)=D_1(k)=D_2(k)=\overline{G}(k)$ be the filling operators,
$\overline{C}(k)$ be the right shift operator.
In this case, from the coupled backward Riccati operator equations \eqref{29}-\eqref{30}, one has that
\begin{equation*}
	\begin{split}
        &P_1(2)=0,\ P_2(2)=0,\\
		&P_1(1)[x_1, x_2,x_3,x_4, \cdots]=-[x_1, x_2,x_3,x_4, \cdots],\\
		&P_2(1)[x_1, x_2,x_3,x_4, \cdots]=[x_1, x_2,x_3,x_4, \cdots],\\
		&P_1(0)[x_1, x_2,x_3,x_4, \cdots]=-\frac{3}{2}[x_1, x_2,x_3,x_4, \cdots]\\
		&\ \ \ \ \ \ \ \ \ \ \ \ \ \ \ \ \ \ \ \ \ \ \ \ \ \ \ \ \ \ \ \ \ \ \ \ \ \ +[\Omega_1x_1, 0,0,0, \cdots],\\
		&P_2(0)[x_1, x_2,x_3,x_4, \cdots]=\frac{3}{2}[x_1, x_2,x_3,x_4, \cdots]\\
		&\ \ \ \ \ \ \ \ \ \ \ \ \ \ \ \ \ \ \ \ \ \ \ \ \ \ \ \ \ \ \ \ \ \ \ \ \ \ +[\Omega_2x_1, 0,0,0, \cdots],\\
		&K_1(1,P_1)=0,\ K_2(1,P_2)=0,\\
		&K_1(0,P_1)[x_1, x_2,x_3,x_4, \cdots]=\Upsilon_1 x_1,\\
		&K_2(0,P_2)[x_1, x_2,x_3,x_4, \cdots]=\Upsilon_2 x_1,
	\end{split}
\end{equation*}
where
$\Omega_1=-\Upsilon_1-2\Upsilon_2-2\Upsilon_1\Upsilon_2-3\Upsilon_2^2,$
$\Omega_2=2\Upsilon_1+\Upsilon_2+2\Upsilon_1\Upsilon_2+(2-\rho^2)\Upsilon_1^2,$
$\Upsilon_1=\frac{1}{3\gamma^2-2}$,
and $\Upsilon_2=-\frac{\gamma^2}{3\gamma^2-2}.$

Take the initial state $x_0=(1, \frac{\sqrt{2}}{2}, (\frac{\sqrt{2}}{2})^2, \cdots)\in \ell^2$.
When $(\gamma,\rho)\in [2,3]\times[0,1]$, by Theorem \ref{theorem-4},
the optimal Nash strategy pair is shown in Fig.\ref{fig3}.
\begin{figure}[htp]
	\centering
	\subfigure[]{
		\includegraphics[height=4cm]{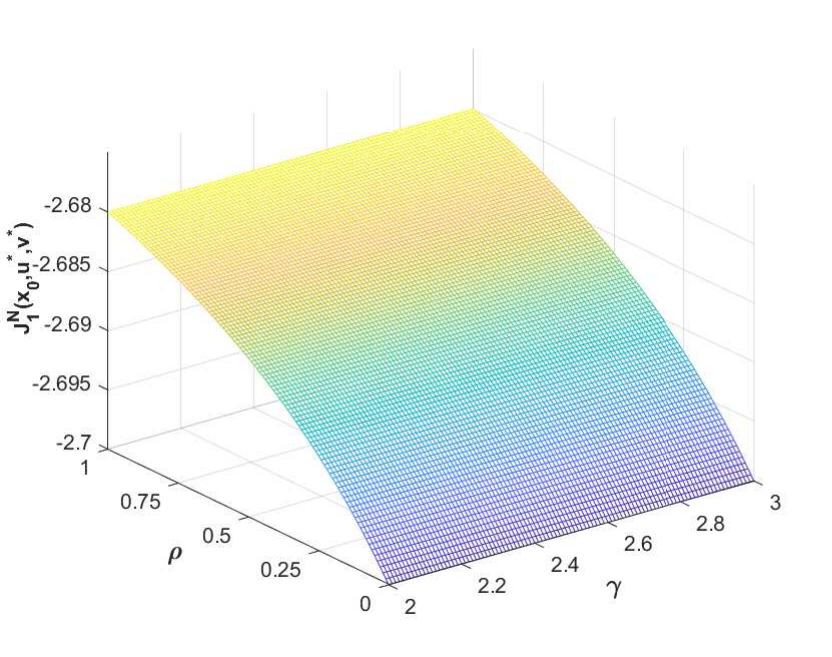}
	}
	\subfigure[]{
		\includegraphics[height=4cm]{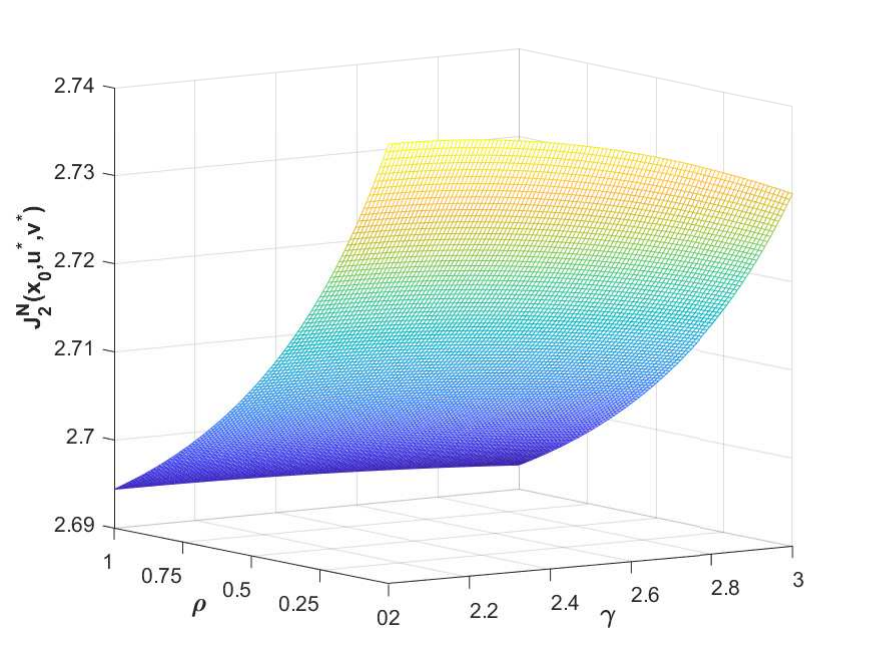}
	}
	\caption{$J^N_1(x_0, u^*, v^*)$ and $J^N_2(x_0, u^*, v^*)$.}
    \label{fig3}	
\end{figure}
\end{exmp}

\section{Conclusion}

In this paper, by leveraging the properties of bounded linear operators and the inner product on Hilbert spaces,
the finite-horizon indefinite LQ-optimal control problem has first been resolved.
Building on this, a sufficient and necessary condition for stochastic bounded real lemma has been given in a real separable Hilbert space.
In addition, the existence of optimal control strategies in a Nash game problem has been discussed,
and a unified design method for the $H_{\infty}$ and $H_2/H_{\infty}$ control has been conceived.
These works construct a theoretical foundation for the finite-horizon indefinite LQ-optimal control and $H_{\infty}$ control theories
in the infinite-dimensional context, and expand the existing knowledge base.
As seen in Examples 1 and 2, in case no settlement can be reached under the finite-dimensional context,
control theories within the framework of Hilbert spaces have fully demonstrated their large application potential in engineering.
Still within the infinite-dimensional setting,
whether the above research results can be extended to the case of infinite-horizon
remains a challenging topic that deserves more effort in our future research.

\section*{Appendix}
\noindent\textbf{PROOF OF THEOREM \ref{theorem-4}.}
Firstly, according to Assumption \ref{assumption-2},
it is noted that via the inner product defined on the suitable Hilbert space,
the performance indices \eqref{25}-\eqref{26} can be rewritten as
    \begin{equation*}
		\begin{split}
			J_1^N(x_0,u,v)=&E\bigg[\sum_{k=0}^{N}\Big(\langle -\overline{C}(k)^*\overline{C}(k)x(k),x(k) \rangle_H\\
			&+\langle \gamma^2 v(k),v(k)\rangle_V-\langle u(k),u(k)\rangle_U\Big)\bigg],\\
			J_2^N(x_0,u,v)=&E\bigg[\sum_{k=0}^{N}\Big(\langle \overline{C}(k)^*\overline{C}(k)x(k),x(k) \rangle_H\\
			&+\langle  u(k),u(k)\rangle_U-\langle \rho^2 v(k),v(k)\rangle_V\Big)\bigg].
		\end{split}
	\end{equation*}

Sufficiency.
Since the coupled backward Riccati operator equations \eqref{29}-\eqref{30}
admit a pair of solutions $(P_1,P_2)\in L^2_{\varepsilon(H)^{{\mathcal{Q}^{+}}}}\times L^2_{\varepsilon(H)^{{\mathcal{Q}^{+}}}}$,
one can set $u^*(k)=K_2(k,P_2)x(k)$ ($k\in \mathcal{Q}$) and substitute it into system \eqref{24} to get
	\begin{equation*}
		\left\{\begin{aligned}
			&x(k+1)=(A(k)+B_2(k)K_2(k,P_2))x(k)\\
			&\ \ \ \ \ \ \ \ +B_1(k)v(k)+((C(k)+D_2(k)K_2(k,P_2))x(k)\\
            &\ \ \ \ \ \ \ \ +D_1(k)v(k))\omega(k),\\
			&z(k)=(\overline{C}(k)+\overline{G}(k)K_2(k,P_2))x(k),\ k\in \mathcal{Q},\\
			&x(0)=x_0\in H.
		\end{aligned}
		\right.
	\end{equation*}
Meanwhile, the performance index \eqref{25} becomes
$J_1^N(x_0,u^*,v)=E[\sum_{k=0}^{N}(\langle -(\overline{C}(k)^*\overline{C}(k)+K_2(k,P_2)^*$
$\cdot K_2(k,P_2))x(k),x(k)\rangle_H+\langle \gamma^2 v(k),v(k)\rangle_V)].$
In \eqref{2}, let $M(k)=-(\overline{C}(k)^*\overline{C}(k)+K_2(k,P_2)^*K_2(k,P_2))$, $L(k)=0$, $R(k)=\gamma^2I_V$ ($k\in \mathcal{Q}$)
and $S(N+1)=0$.
Because $\{P_1(k),\ k\in \mathcal{Q}^{+}\}$ satisfies the coupled backward Riccati operator equation \eqref{29},
it can be carried out by Theorem \ref{theorem-1} that for any $v\in L^2_{v^{\mathcal{Q}}}$,
	\begin{equation}\label{25pl}
		\begin{aligned}
			J_1^N(x_0,u^*,v)\geq J_1^N(x_0,u^*,v^*)=\langle P_1(0)x_0,x_0 \rangle_H,
		\end{aligned}
	\end{equation}
where $v^*=\{v^*(k)=K_1(k,P_1)x(k),\ k\in \mathcal{Q}\}$.
	
Similarly, substituting $v^*(k)=K_1(k,P_1)x(k)$ ($k\in \mathcal{Q}$) into system \eqref{24}
and the parameterized quadratic performance index $J_2^N(x_0,u,v^*)$, we have that
	\begin{equation*}
		\left\{\begin{aligned}
			&x(k+1)=(A(k)+B_1(k)K_1(k,P_1))x(k)\\
			&\ \ \ \ \ \ \ \ +B_2(k)u(k)+((C(k)+D_1(k)K_1(k,P_1))x(k)\\
			&\ \ \ \ \ \ \ \ +D_2(k)u(k))\omega(k),\\
			&z(k)=\overline{C}(k)x(k)+\overline{G}(k)u(k),\ k\in \mathcal{Q},\\
			&x(0)=x_0\in H,
		\end{aligned}
		\right.
	\end{equation*}
and
$J_2^N(x_0,u,v^*)=E[\sum_{k=0}^{N}(\langle [\overline{C}(k)^*\overline{C}(k)-\rho^2K_1(k,P_1)^*$
$\cdot K_1(k,P_1)]x(k), x(k) \rangle_H+\langle u(k),u(k)\rangle_U)].$
Subsequently, applying Theorem \ref{theorem-1}, one can deduce that
	\begin{equation}\label{36}
		\begin{aligned}
			J_2^N(x_0,u,v^*)\geq J_2^N(x_0,u^*,v^*)=\langle P_2(0)x_0,x_0 \rangle_H,
		\end{aligned}
	\end{equation}
where $u^*=\{u^*(k)=K_2(k,P_2)x(k),\ k\in \mathcal{Q}\}$.
By combining \eqref{25pl} with \eqref{36},
we draw a conclusion that $(u^*,v^*)$ represents the optimal Nash strategy pair correlating with $J_1^N(x_0,u,v)$ and $J_2^N(x_0,u,v)$.
	
Necessity.
Assume that the Nash game problem \eqref{27}-\eqref{28} associated with system \eqref{24}
is solved via a linear feedback Nash equilibrium $(u^*,v^*)$.
Without loss of generality, for any $k\in \mathcal{Q}$, denote
$u^*(k)=\overline{K}_2(k)x(k)$, $v^*(k)=\overline{K}_1(k)x(k)$,
where $x(k)$ is the state trajectory of system \eqref{24},
$\overline{K}_2(k)$ and $\overline{K}_1(k)$ are linear feedback gains.
Therefore, according to Definition \ref{definition-5}, $v^*$ solves the following LQ-optimal problem:
	\begin{equation*}
		\begin{split}
			\inf&_{v \in L^2_{V^{\mathcal{Q}}}}\bigg\{J_1^N(x_0,u^*,v)=E\bigg[\sum_{k=0}^{N}\Big(\langle -(\overline{C}(k)^*\overline{C}(k)
            \\
			&+\overline{K}_2(k)^*\overline{K}_2(k))x(k),x(k) \rangle_H+\langle \gamma^2 v(k),v(k)\rangle_V\Big)\bigg]\bigg\}
		\end{split}
	\end{equation*}
subject to
	\begin{equation}\label{37}
		\left\{\begin{aligned}
			&x(k+1)=(A(k)+B_2(k)\overline{K}_2(k))x(k)+B_1(k)v(k)\\
			&\ \ \ +((C(k)+D_2(k)\overline{K}_2(k))x(k)+D_1(k)v(k))\omega(k),\\
			&z(k)=(\overline{C}(k)+\overline{G}(k)\overline{K}_2(k))x(k),\ k\in \mathcal{Q},\\
			&x(0)=x_0\in H.
		\end{aligned}
		\right.
	\end{equation}	
For system \eqref{37}, set $x_0=0$.
$J_1^N(0,u^*,v^*)=0$ is then derived from $v^*=0$.
Thereby, for any $v \in L^2_{V^{\mathcal{Q}}},\ v\neq 0$, we have that $J_1^N(x_0,u^*,v)>0$.
Hence, by Proposition \ref{proposition-2} and Theorem \ref{theorem-3},
it yields that the following backward operator equation:
	\begin{equation}\label{38}
		\left\{\begin{aligned}
			&\overline{P}_1(k)=[A(k)+B_2(k)\overline{K}_2(k)]^*\overline{P}_1(k+1)\\
			&\ \ \cdot [A(k)+B_2(k)\overline{K}_2(k)]+[C(k)+D_2(k)\overline{K}_2(k)]^*\\
			&\ \ \cdot\overline{P}_1(k+1)[C(k)+D_2(k)\overline{K}_2(k)] \\
			&\ \ -\overline{K}_2(k)^*\overline{K}_2(k)-\overline{C}(k)^*\overline{C}(k)\\
			&\ \ -\mathcal{G}_1(k,\overline{P}_1)^*\mathcal{R}_1(k,\overline{P}_1)^{-1}
             \mathcal{G}_1(k,\overline{P}_1),\\
            &\mathcal{R}_1(k,\overline{P}_1)^{-1}>0,\ k\in\mathcal{Q},\\
			&\overline{P}_1(N+1)=0
		\end{aligned}
		\right.
	\end{equation}
admits a solution $\{\overline{P}_1(k),\ k\in\mathcal{Q}^{+}\}\in L^2_{\varepsilon(H)^{{\mathcal{Q}^{+}}}}$.
Further, by applying Theorem \ref{theorem-1}, the unique optimal linear feedback control is
$v^*(k)=-\mathcal{R}_1(k,\overline{P}_1)^{-1}\mathcal{G}_1(k,\overline{P}_1)x(k)$,
i.e. $\overline{K}_1(k)=K_1(k,\overline{P}_1)$.
	
On the other side, by Definition \ref{definition-5} again,
one confirms that the following LQ-optimal problem:
	\begin{equation*}
		\begin{split}
			\inf&_{u\in L^2_{U^{\mathcal{Q}}}}\bigg\{J_2^N(x_0,u,v^*)=E\bigg[\sum_{k=0}^{N}\Big(\langle (\overline{C}(k)^*\overline{C}(k)\\
                                &-\rho^2K_1(k,\overline{P}_1)^*K_1(k,\overline{P}_1))x(k),x(k) \rangle_H\\
                                &+\langle u(k),u(k)\rangle_U\Big)\bigg]\bigg\}
		\end{split}
	\end{equation*}
subject to
	\begin{equation*}
		\left\{\begin{aligned}
			&x(k+1)=(A(k)+B_1(k)K_1(k,\overline{P}_1))x(k)+B_2(k)u(k)\\
			&\ \ \ \ \ \ \ \ \ \ \ \ \ \ \ \ +((C(k)+D_1(k)K_1(k,\overline{P}_1))x(k)\\
			&\ \ \ \ \ \ \ \ \ \ \ \ \ \ \ \ +D_2(k)u(k))\omega(k),\\
			&z(k)=\overline{C}(k)x(k)+\overline{G}(k)u(k),\ k\in \mathcal{Q},\\
			&x(0)=x_0\in H
		\end{aligned}
		\right.
	\end{equation*}
is solved by the linear state feedback control law $\{u^*(k)=\overline{K}_2(k)x(k),\ k\in\mathcal{Q}\}$.
In \eqref{2}, set $M(k)=\overline{C}(k)^*\overline{C}(k)-\rho^2K_1(k,\overline{P}_1)^*K_1(k,\overline{P}_1)$,
$L(k)=0$, $R(k)=I_U$ ($k\in \mathcal{Q}$) and $S(N+1)=0$.
From Theorem \ref{theorem-1}, it reveals that the following backward operator equation:
\begin{equation}\label{39}
	\left\{\begin{aligned}
		&\overline{P}_2(k)=[A(k)+B_1(k)K_1(k,\overline{P}_1)]^*\overline{P}_2(k+1)\\
		&\ \ \cdot [A(k)+B_1(k)K_1(k,\overline{P}_1)]+[C(k)+D_1(k)\\
		&\ \ \cdot K_1(k,\overline{P}_1)]^*\overline{P}_2(k+1)[C(k)+D_1(k)K_1(k,\overline{P}_1)]\\
		&\ \ -\rho^2K_1(k,\overline{P}_1)^*K_1(k,\overline{P}_1)+\overline{C}(k)^*\overline{C}(k)\\
		&\ \ -\mathcal{G}_2(k,\overline{P}_2)^*\mathcal{R}_2(k,\overline{P}_2)^{-1}\mathcal{G}_2(k,\overline{P}_2),\\
        &\mathcal{R}_2(k,\overline{P}_2)^{-1}>0,\ k\in\mathcal{Q},\\
		&\overline{P}_2(N+1)=0
	\end{aligned}
	\right.
\end{equation}
has a solution $\{\overline{P}_2(k),\ k\in\mathcal{Q}^{+}\}\in L^2_{\varepsilon(H)^{{\mathcal{Q}^{+}}}}$
and $u^*(k)=-\mathcal{R}_2(k,\overline{P}_2)^{-1}\mathcal{G}_2(k,\overline{P}_2)x(k)$ is the unique optimal control,
which implies that $\overline{K}_2(k)=K_2(k,\overline{P}_2)$.
Incorporating this result with before,
it demonstrates that \eqref{38} and \eqref{39} are exactly the same with \eqref{29} and \eqref{30},
also proves that $K_1(k,\overline{P}_1)=K_1(k,P_1)$ and $K_2(k,\overline{P}_2)=K_2(k,P_2)$.
Thus, the proof of Theorem \ref{theorem-4} is all completed.

%\begin{ack}                               % Place acknowledgements
%\end{ack}

%\begin{figure}
%\begin{center}
%\includegraphics[height=4cm]{jcaesar.eps}    % The printed column
%\caption{Gaius Julius Caesar, 100--44 B.C.}  % width is 8.4 cm.
%\label{fig1}                                 % Size the figures
%\end{center}                                 % accordingly.
%\end{figure}

% OR

%\begin{figure}
%\begin{center}
%\epsfig{file=jcaesar,width=7cm}
%\caption{Gaius Julius Caesar, 100--44 B.C.}
%\label{fig1}
%\end{center}
%\end{figure}

%\bibliographystyle{plainnat}        % Include this if you use bibtex
%\bibliography{autosam}           % and a bib file to produce the
                                 % bibliography (preferred). The
                                 % correct style is generated by
                                 % Elsevier at the time of printing.

\end{document}